\documentclass[leqno,11pt]{article}
\usepackage{amssymb,amsmath,mathabx}

\usepackage{pgf,ifpdf,graphics,xcolor}

\usepackage{amsfonts,amssymb,euscript,epsfig,color}
\usepackage[all]{xy}

\newcounter{daggerfootnote}

\newtheorem{theo}{Theorem}[section]
\newtheorem{defi}[theo]{Definition}

\newtheorem{lem}[theo]{Lemma}

\newtheorem{prop}[theo]{Proposition}

\newtheorem{rem}[theo]{Remark}
\newtheorem{coro}[theo]{Corollary}

\newtheorem{exam}[theo]{Example}

\newcommand{\bgot}{\ensuremath{\mathfrak{b}}}

\newcommand{\ugot}{\ensuremath{\mathfrak{u}}}
\newcommand{\kgot}{\ensuremath{\mathfrak{k}}}

\newcommand{\ggot}{\ensuremath{\mathfrak{g}}}
\newcommand{\ngot}{\ensuremath{\mathfrak{n}}}

\newcommand{\tgot}{\ensuremath{\mathfrak{t}}}

\newcommand{\Rgot}{\ensuremath{\mathfrak{R}}}
\newcommand{\Sgot}{\ensuremath{\mathfrak{S}}}
\newcommand{\Xgot}{\ensuremath{\mathfrak{X}}}

\newcommand{\Acal}{\ensuremath{\mathcal{A}}}

\newcommand{\Ccal}{\ensuremath{\mathcal{C}}}

\newcommand{\Ecal}{\ensuremath{\mathcal{E}}}
\newcommand{\Fcal}{\ensuremath{\mathcal{F}}}

\newcommand{\Lcal}{\ensuremath{\mathcal{L}}}
\newcommand{\Mcal}{\ensuremath{\mathcal{M}}}

\newcommand{\Ocal}{\ensuremath{\mathcal{O}}}
\newcommand{\Pcal}{\ensuremath{\mathcal{P}}}

\newcommand{\Vcal}{\ensuremath{\mathcal{V}}}

\newcommand{\Xcal}{\ensuremath{\mathcal{X}}}

\newcommand{\Ebb}{\ensuremath{\mathbb{E}}}

\newcommand{\Qbb}{\ensuremath{\mathbb{Q}}}
\newcommand{\Nbb}{\ensuremath{\mathbb{N}}}
\newcommand{\Pbb}{\ensuremath{\mathbb{P}}}
\newcommand{\Cbb}{\ensuremath{\mathbb{C}}}
\newcommand{\Gbb}{\ensuremath{\mathbb{G}}}
\newcommand{\Rbb}{\ensuremath{\mathbb{R}}}
\newcommand{\Zbb}{\ensuremath{\mathbb{Z}}}

\newcommand{\T}{\ensuremath{\hbox{\bf T}}}

\newcommand{\End}{\ensuremath{\hbox{\rm End}}}

\newcommand{\tr}{\operatorname{Tr}}

\newcommand{\Eul}{\operatorname{Eul}}

\newcommand{\horn}{\ensuremath{\hbox{\rm Horn}}}
\newcommand{\Horn}{\ensuremath{{\rm Horn}_{\rm hol}}}

\newcommand{\maxx}{\operatorname{\hbox{\rm \tiny max}}}

\def \tI {\widetilde{I}}

\begin{document}

\title{$\horn(p,q)$}

\author{Paul-Emile Paradan\footnote{IMAG, Univ Montpellier, CNRS, email : paul-emile.paradan@umontpellier.fr}}

\maketitle

\date{}


\begin{abstract}
In this article, we obtain a recursive description of the Horn cone $\horn(p,q)$ with respect to the integers $p$ and $q$, as in the classical Horn's conjecture.
\end{abstract}


\tableofcontents

\section{Introduction}

When $G$ is Lie group, a natural problem is to understand how the sum of two adjoint orbits decomposes into a union of adjoint orbits. 
Let $\ggot$ be the the Lie algebra of $G$ and let $\ggot/G$ be the set of adjoint orbits. The {\em Horn cone}  is defined as follows
$$
\horn(G)=\{(\Ocal,\Ocal',\Ocal'')\in (\ggot/G)^3,\ \Ocal''\subset \Ocal+\Ocal'\}.
$$

Consider the case where $G$ is a compact connected Lie group. Let $T \subset G$ be a maximal torus with Lie algebra $\tgot$. 
The set $\ggot/G$ admits a canonical identification with a Weyl chamber $\tgot_{\geq 0}\subset\tgot$. In this setting, 
the Horn cone $\horn(G)\subset (\tgot_{\geq 0})^3$ has been at the center of numerous studies 
\cite{Horn,Klyachko,Knutson-Tao-99,Berenstein-Sjamaar,Belkale06,BK1,Ka-Le-Mi,Ressayre10} 
that we summarize by the following theorem. We refer the reader to the survey articles \cite{Brion-Bourbaki-11,Kumar14} for details.

\begin{theo}
If $G$ is a compact connected Lie group, $\horn(G)$ is a polyhedral convex cone and one can parametrize the equation of its facets by cohomological means.
\end{theo}

\subsection{Horn's conjecture}

When $G$ is the unitary group $U(n)$,  the convex polyhedral cone\footnote{We note $\horn(U(n))$ simply by $\horn(n)$.} $\horn(n)$ has a nice feature which was predicted by A. Horn in the 60s : it admits a recursive description relative to the integer $n\geq 1$ \cite{Horn}. 

Denote the set of cardinality $r$ subsets $I=\{i_1<i_2<\cdots<i_r\}$ of $[n]=\{1,\ldots,n\}$ by $\Pcal^n_r$. To each $I\in \Pcal^n_r$
we associate a weakly decreasing sequence of non-negative integers 
\begin{equation}\label{eq:lambda-I}
\lambda(I)=(\lambda_1\geq\lambda_2\cdots\geq \lambda_r)\ \in\ \Zbb_{\geq 0}^r
\end{equation}
where $\lambda_a= n-r+a-i_a$ for $a\in [r]$.

Let ${\bf d}:\Rbb^n\to \ugot(n)$ be the map that sends $X=(x_1,\ldots,x_n)$ to the diagonal matrix ${\bf d}_X={\rm Diag}(ix_1,\ldots,ix_n)$. 
The map ${\bf d}$ induces a one to one correspondence between $\Ccal_n=\{(x_1\geq\cdots\geq x_n)\}\subset\Rbb^n$ and the set of $U(n)$-adjoint orbits.
If $X=(x_1,\ldots,x_n)\in\Rbb^n$ and $I\subset [n]$, we define $|\,X\,|_I=\sum_{i\in I}x_i$ and $|\,X\,|=\sum_{i=1}^n x_i$.

\begin{defi} \label{defi:horn-cone} Let $n\geq 1$.
$$
\horn(n)=\{(A,B,C)\in(\Ccal_n)^3,\ U(n){\bf d}_C\subset U(n){\bf d}_A+U(n){\bf d}_B\}.
$$
\end{defi}

\medskip

The following Horn's conjecture \cite{Horn} was settled in the affirmative by combining the work of A. Klyachko \cite{Klyachko} with the work of A. Knutson and T. Tao 
\cite{Knutson-Tao-99} on the ``saturation'' problem. We refer the reader to Fulton's survey article \cite{Fulton-00} for details.

\begin{theo}[Horn's conjecture]

An element $(A,B,C)\in(\Ccal_n)^3$ belongs to $\horn(n)$ if and only if the following conditions holds
\begin{itemize}
\item $|\,A\,|+|\,B\,|=|\,C\,|$, 
\item $\forall r\in [n-1]$, $\forall I,J,K\in\Pcal^n_r$, we have 
$$
|\,A\,|_{I}+|\,B\,|_{J} \leq |\,C\,|_{K}\quad {\rm if}\quad (\lambda(I),\lambda(J),\lambda(K))\in\horn(r).
$$
\end{itemize}
\end{theo}

\subsection{Holomorphic Horn cone $\Horn(p,q)$}

Let $p\geq q\geq 1$. We begin by recalling the definition of the holomorphic Horn cone $\Horn(p,q)$ associated with the pseudo-unitary group $U(p,q)$.

The Lie group $U(p,q)\subset GL_{p+q}(\Cbb)$ is defined by the relations  
$g{\rm Id}_{p,q}g^*={\rm Id}_{p,q}$, where ${\rm Id}_{p,q}$ is the diagonal matrice ${\rm Diag}({\rm Id}_p,-{\rm Id}_q)$. 
The Lie algebra $\ugot(p,q)$ of $U(p,q)$ admits the following invariant convex cone 
$$
\Ccal(p,q)=\left\{X\in\ugot(p,q),\ {\rm Im}(\tr( gXg^{-1} {\rm Id}_{p,q})) \geq 0,\ \forall g\in U(p,q)\right\}.
$$
Let us consider
$$
\Ccal_{p,q}=\left\{x\in \Rbb^p\times\Rbb^q, x_1\geq \cdots \geq x_p > x_{p+1}\geq \cdots \geq x_{p+q}\right\}\subset\Ccal_p\times\Ccal_q
$$
and the map ${\bf d}: \Rbb^p\times\Rbb^q\to\ugot(p,q)$. A well-know result says that for any $U(p,q)$-orbit $\Ocal$ contained in the interior of $\Ccal(p,q)$, 
there exists a unique $X\in \Ccal_{p,q}$ such that $\Ocal=U(p,q){\bf d}_X$ (see \cite{Vinberg80,Paneitz83}). In other words, the map ${\bf d}$ realizes a one to one map between $\Ccal_{p,q}$ and the set of $U(p,q)$-orbits in 
the interior of the invariant convex cone $\Ccal(p,q)$. The holomorphic Horn cone is then defined as follows : 
$$
\Horn(p,q)=\left\{(A,B,C)\in(\Ccal_{p,q})^3,\ U(p,q){\bf d}_C\subset U(p,q){\bf d}_A+U(p,q){\bf d}_B\right\}.
$$

In a companion paper \cite{pep:hol-horn}, we have proved that  $\Horn(p,q)$ is a closed convex cone of $(\Ccal_{p,q})^3$, and we have explained 
a way to compute it. In order to detail this result, we need some additional notations. For any $n\geq 1$, we consider the semigroup
$\wedge_n^+=\{(\lambda_1\geq \cdots\geq\lambda_n)\}\subset\Zbb^n$. If $\lambda=(\lambda',\lambda'')\in \wedge_p^+\times\wedge_q^+$, then
$V_\lambda:=V^{U(p)}_{\lambda'}\otimes V^{U(q)}_{\lambda''}$ denotes the irreducible representation of $U(p)\times U(q)$ with highest weight $\lambda$. 
We denote by $M_{p,q}$ the vector space of $p\times q$ complex matrices, and by ${\rm Sym}(M_{p,q})$ the symmetric algebra of $M_{p,q}$.

If $H$ is a representation of $U(p)\times U(q)$, we denote by $[V_\nu : H]$ the multiplicity of $V_\nu$ in $H$.

\begin{defi}\label{defi:horn-p-q} 
\begin{enumerate}
\item $\horn^{\Zbb}(p,q)$ is the semigroup of $(\wedge_p^+\times\wedge_q^+)^3$ defined by the conditions:
$$
 (\lambda,\mu,\nu)\in \horn^{\Zbb}(p,q)\Longleftrightarrow \left[V_\nu\,:\,V_\lambda\otimes V_\mu\otimes {\rm Sym}(M_{p,q})\right]\neq 0.
$$
\item $\horn(p,q)$ is the convex cone of $(\Ccal_p\times\Ccal_q)^3$ defined as the closure of 
$\Qbb^{>0}\cdot \horn^{\Zbb}(p,q)$.
\end{enumerate}
\end{defi}

The following result is proved in \cite{pep:hol-horn}.

\begin{theo}
We have 
$$
\Horn(p,q)=\horn(p,q)\,\bigcap \,(\Ccal_{p,q})^3.
$$
\end{theo}

\subsection{Statement of the main result}

We now explain the main purpose of this paper that concerns a recursive description of the convex polyhedral cones $\horn(p,q)$ as in Horn's conjecture.
We need another notations.

\begin{enumerate}
\item If $A=(A',A'')\in\Rbb^p\times \Rbb^q$ and $I=I'\times I''\subset [p] \times [q]$, we define
$|\,A\,|_I=|\,A'\,|_{I'}+|\,A''\,|_{I''}$ and  $|\,A\,|=|\,A'\,|+|\,A''\,|$.
\item If $I=I'\times I''\subset [n]\times [m]$ then $\lambda(I)=(\lambda(I'),\lambda(I''))\in \wedge_n^+\times \wedge_m^+$.
\item Let $\mathbf{1}_n=(1,\ldots,1)\in\Zbb^n$.
\end{enumerate}

The main result of this paper is the following theorem.

\begin{theo} \label{theo:main} Let $p\geq q\geq 1$. 
An element $(A,B,C)\in (\Ccal_p\times\Ccal_q)^3$ belongs to 
$\horn(p,q)$ if and only if the following conditions holds:
\begin{itemize}
\item $\boxed{|\,A\,| + |\,B\,| = |\,C\,|}$.
\item $\boxed{|\,A'\,| + |\,B'\,| \leq |\,C'\,|}$.
\item For any $r\in [p\!-\!1]$, for any $I',J',K'\in \Pcal^p_r$, we have :
\begin{align*}
\boxed{|\,A'\,|_{I'} + |\,B'\,|_{J'} \leq |\,C'\,|_{K'}} \quad &{\rm if} \ (\lambda(I'),\lambda(J'),\lambda(K'))\in\horn(r).\\
\boxed{|\,A'\,|_{I'} + |\,B'\,|_{J'} \geq |\,C'\,|_{K'}} \quad & {\rm if} \
(\lambda(I'),\lambda(J'),\lambda(K')+(q\!+\!p\!-\!r)\mathbf{1}_r)
\in\horn(r).
\end{align*}
\item For any $s\in [q\!-\!1]$, for any $I'',J'',K''\in \Pcal^q_s$, we have :
\begin{align*}
\boxed{|\,A''\,|_{I''} + |\,B''\,|_{J''} \geq |\,C''\,|_{K''}} \quad &{\rm if} \ (\lambda(I''),\lambda(J''),\lambda(K'')+(q\!-\!s)\mathbf{1}_s)\in\horn(s).\\
\boxed{|\,A''\,|_{I''} + |\,B''\,|_{J''} \leq |\,C''\,|_{K''}} \quad & {\rm if} \ 
(\lambda(I''),\lambda(J''),\lambda(K'')-p\mathbf{1}_s)
\in\horn(s).
\end{align*}
\item For any $(r,s)\in [p\!-\!1]\times [q\!-\!1]$ with $r\geq s$, for any $I,J,K\in \Pcal^p_r\times\Pcal^q_s$, we have 
$$
\boxed{|\,A\,|_{I} + |\,B\,|_{J}  \leq |\,C\,|_{K}} \quad {\rm if} \ 
\big(\lambda(I),\lambda(J),\lambda(K)+ (0,(r-p)\mathbf{1}_s)\big)\in \horn(r,s).
$$
\end{itemize}
\end{theo}

\medskip

\subsection{Examples}

The convex cones $\horn(1,1)$, $\horn(2,1)$ and $\horn(2,2)$ admit the following descriptions. 

\begin{prop}
An element $(A,B,C)\in (\Rbb\times\Rbb)^3$ belongs to 
$\horn(1,1)$ if and only if the following conditions holds:
$$
\boxed{a_1+a_2 + b_1 + b_2 = c_1+c_2}
$$
$$
\boxed{a_1+b_1\leq c_1}
$$
\end{prop}

\medskip

\begin{prop}
An element $(A,B,C)\in (\Ccal_2\times\Rbb)^3$ belongs to 
$\horn(2,1)$ if and only if  the following conditions holds:
$$
\boxed{a_1 + a_2 + a_3 +b_1+ b_2 + b_3  = c_1 + c_2  +c_3}
$$
$$
\boxed{a_1 + a_2  +b_1+ b_2   \leq c_1 + c_2  }
$$
\begin{equation*}
\boxed{
\begin{array}{rcl}
a_2+b_2 &\leq & c_2\\
a_2+b_1 &\leq & c_1\\
a_1+b_2 &\leq & c_1\\
a_1+b_1 &\geq & c_2
\end{array}
}
\end{equation*}

\end{prop}

\medskip

\begin{prop}
An element $(A,B,C)\in (\Ccal_2\times\Ccal_2)^3$ belongs to 
$\horn(2,2)$ if and only if the following conditions holds:

$$
\boxed{a_1+a_2+a_3+a_4+b_1+b_2+b_3+b_4=c_1+c_2+c_3+c_4}
$$
$$
\boxed{a_1+a_2+b_1+b_2\leq c_1+c_2}
$$
\begin{equation*}
\boxed{
\begin{array}{rcl}
a_2+b_2&\leq & c_2\\
a_2+b_1&\leq & c_1\\
a_1+b_2&\leq & c_1
\end{array}
}
\end{equation*}
\begin{equation*}
\boxed{
\begin{array}{rcl}
a_3+b_3 &\geq & c_3\\
a_3+b_4 &\geq & c_4\\
a_4+b_3 &\geq & c_4
\end{array}
}
\end{equation*}
\begin{equation*}
\boxed{
\begin{array}{rcl}
a_2+a_4+b_2+b_4 &\leq & c_1+c_4 \\
a_2+a_4+b_2+b_4 &\leq & c_2+c_3\\
a_2+a_4+b_1+b_4 &\leq & c_1+c_3\\
a_1+a_4+b_2+b_4 &\leq & c_1+c_3\\
a_2+a_4+b_2+b_3 &\leq & c_1+c_3\\
a_2+a_3+b_2+b_4 &\leq & c_1+c_3
\end{array}
}
\end{equation*}
\end{prop}

\subsection{Outline of the article}

The recursive description of $\horn(p,q)$ is obtained by studying the Hamiltonian action of $(U(p)\times U(q))^3$ on the manifold\footnote{We use the notation $GL_n$ for 
the Lie group $GL(\Cbb^n)$.}
$(GL_p \times GL_q)^2 \times \Cbb^p\otimes \Cbb^q$. Let ${\rm S}(p,q)\subset
(\Ccal_p\times \Ccal_q)^3$ be the corresponding Kirwan polyhedron.

In \S \ref{sec:K-trois-E}, we study the general framework of a Hamiltonian action of a compact Lie group $K^3$ on 
$(K_\Cbb\times K_\Cbb)^2\times E$ : here $E$ is a $K$-module such that  the coordinate ring $\Cbb[E]$ does 
not admit non-constant invariant vectors. We explain how to parameterize the facets of the Kirwan polyhedron 
$\Delta((K_\Cbb\times K_\Cbb)^2\times E)$ in terms of Ressayre's data \cite{pep-ressayre-pair}. This parametrization requires two steps~: 
determination of the admissible elements which are the potential vectors orthogonal to the facets, and computation of cohomological conditions on  flag varieties.

In \S \ref{sec:saturated semi-group}, we check that the semigroup $\horn^\Zbb(p,q)$ is saturated. It is a direct consequence of the 
Darksen-Weyman saturation theorem \cite{Derksen-Weyman}.

In \S \ref{sec:polyhedron-s-p-q}, we determine the admissible elements relative to the action of $(U(p)\times U(q))^3$ on
$(GL_p \times GL_q)^2 \times \Cbb^p\otimes \Cbb^q$, and we detailed the cohomological conditions in this particular case. 
The formulas need the computation of certain Euler classes which we carry over to \S \ref{sec:euler-class}.

In \S \ref{sec:s-p-q-facets}, we calculate (recursively) the facets of the Kirwan polyhedron ${\rm S}(p,q)$. In the last subsection, we complete the proof of our main result.

\medskip

\subsection*{Acknowledgements} I wish to thank Mich\`ele Vergne for our discussions on this subject and for pointing to my attention the Derksen-Weyman saturation theorem.

\medskip

\section{The $K^3$-manifold $K_\Cbb\times K_\Cbb\times E$}\label{sec:K-trois-E}

In this section, we briefly recall the result of \S 6 of \cite{pep-ressayre-pair} concerning the parametrization of the facets of Kirwan polyhedrons in 
terms of Ressayre's data. 

Let $K$ be a compact connected Lie group with complexification $K_\Cbb$. Let $T\subset K$ be a maximal torus with Lie algebra $\tgot$. We consider the lattice $\wedge:=\frac{1}{2\pi}\ker(\exp:\tgot\to T)$ and the dual lattice $\wedge^*\subset \tgot^*$ defined by $\wedge^*=\hom(\wedge,\Zbb)$. 
We remark that $i\eta$ is a differential of a character of $T$ if and only if $\eta\in\wedge^*$. 
The $\Qbb$-vector space generated by the lattice $\wedge^*$ is denoted by $\tgot^*_\Qbb$: the vectors belonging to 
$\tgot^*_\Qbb$ are designed as rational. Let $\tgot^*_{\geq 0}$ be a Weyl chamber. The set $\wedge^*_+:=\wedge^*\cap\tgot^*_{\geq 0}$ parametrizes the  irreducible representations of $K$: for any $\mu\in\wedge^*_+$, we denote by $V_\mu$ the irreducible representation of $K$ with highest weight $\mu$. 

When $K$ acts linearly on a vector space $H$, we denote by $H^K$ the subspace of invariant vectors under the $K$-action.

Let $E$ be a $K$-module such that $\Cbb[E]^K=\Cbb$ : hence the coordinate ring $\Cbb[E]$ has finite $K$-multiplicities. We consider the following $K\times K\times K$ action on the affine variety  $K_\Cbb\times K_\Cbb\times E$ : 
$$
(k_1,k_2,k_3)\cdot (x,y,v)= (k_1 x k_3^{-1}, k_2 y k_3^{-1},k_3 v).
$$
The  coordinate ring $\Cbb[K_\Cbb\times K_\Cbb\times E]$, viewed as a $K^3$-module, admits the following decomposition 
$$
\Cbb[K_\Cbb\times K_\Cbb\times E]=\sum_{\lambda,\mu,\nu\in\wedge^*_+}m_E(\lambda,\mu,\nu) \ V^{1}_\lambda \otimes V^{2}_\mu \otimes V^{3}_\nu,
$$
where $m_E(\lambda,\mu,\nu)  =\dim [V_\lambda\otimes V_\mu\otimes V_\nu\otimes {\rm Sym}(E)]^{K}$.

\begin{defi}We define the following sets :
\begin{itemize}
\item The semigroup $\Delta^\Zbb(K_\Cbb\times K_\Cbb\times E)\subset (\wedge^*_+)^3$ is defined as follows: 
$(\lambda,\mu,\nu)\in \Delta^\Zbb(K_\Cbb\times K_\Cbb\times E)\Longleftrightarrow m_E(\lambda,\mu,\nu)\neq 0$.
\item The convex cone $\Delta(K_\Cbb\times K_\Cbb\times E)\subset(\tgot^*_{\geq 0})^3$ is  the closure of \break 
$\Qbb^{>0}\cdot\Delta^\Zbb(K_\Cbb\times K_\Cbb\times E)$.
\end{itemize}
\end{defi}

Let us explain why the complex $K^3$-manifold $N=K_\Cbb\times K_\Cbb\times E$ admits a symplectic structure $\Omega_N$ compatible with 
the complex structure, and a moment map $\Phi:N\to(\kgot^*)^3$ associated to the action of $K^3$ on $(N,\Omega_N)$.

Let $h_E$ be a $K$-invariant hermitian structure on $E$. We equip $E$ with the $2$-form 
$\Omega_E=-{\rm Im}(h_E)$. The moment map $\Phi_E$ relative to the action of $K$ on the symplectic vector space $(E,\Omega_E)$ is defined by 
$$
\langle\Phi_E(v),X\rangle=\tfrac{1}{2}\Omega_E(Xv,v),\quad \forall v\in V,\ \forall X\in\kgot.
$$
The hypothesis $\Cbb[E]^K=\Cbb$ implies that $\Phi_E$ is a proper map.

There is a diffeomorphism of the cotangent bundle $\T^*K$ with $K_\Cbb$ defined as follows. 
We identify $\T^*K$ with $K\times \kgot^*$ by means of left-translation and then with $K\times \kgot$ by means of an invariant inner product on $\kgot$. 
The map $\varphi:K\times \kgot \to K_\Cbb$ given by $\varphi(k,X)=ke^{iX}$ is a diffeomorphism. If we use $\varphi$ to transport the canonical symplectic $2$-form of   
$\T^*K$ to $K_\Cbb$, then the resulting $2$-form $\Omega_{K_\Cbb}$ on $K_\Cbb$ is compatible with the complex structure (see \cite{Hall87}, \S 3).

Finally, the $K^3$-manifold $K_\Cbb\times K_\Cbb\times E\simeq \T^* K\times \T^* K\times E$ carries the symplectic 2-form 
$\Omega_N:=\Omega_{K_\Cbb}\times\Omega_{K_\Cbb}\times \Omega_E$ which is compatible with the complex structure. The moment map relative to the $K^3$-action 
on $(N,\Omega_N)$ is the proper map $\Phi=\Phi_1\oplus\Phi_2\oplus\Phi_3 : \T^* K\times \T^* K\times E \to \kgot^*\oplus \kgot^*\oplus \kgot^*$ defined by
\begin{equation}\label{eq:momentcotangent}
\Phi(g_1,\xi_1,g_2,\xi_2, v)=(-g_1\xi_1,-g_2\xi_2,\xi_1+\xi_2+\Phi_E(v)).
\end{equation}

By definition, the Kirwan polyhedron $\Delta(\T^* K\times \T^* K\times E)$ is the intersection of the image of $\Phi$ with $(\tgot^*_{\geq 0})^3$. 
The following result is classical (see Theorem 4.9 in \cite{Sjamaar-98}).
\begin{prop}\label{prop:delta=delta}
The Kirwan polyhedron $\Delta(\T^* K\times \T^* K\times E)$ is equal to $\Delta(K_\Cbb\times K_\Cbb\times E)$.
\end{prop}

\subsection{Admissible elements}

\begin{defi} 
When a Lie group $G$ acts on a manifold $N$, the stabilizer subgroup of $n\in N$ is denoted by $G_n=\{g\in G,gn=n\}$, and its Lie algebra by $\ggot_n$.
Let us define $\dim_G(\Xcal)=\min_{n\in \Xcal} \dim(\ggot_n)$ for any subset $\Xcal\subset N$. 
\end{defi}

We start by introducing the notion of admissible elements. The group $\hom(U(1),T)$ admits a natural identification with the lattice 
$\wedge:=\frac{1}{2\pi}\ker(\exp:\tgot\to T)$. A vector $\gamma\in\tgot$ is called rational if it belongs to the $\Qbb$-vector space 
$\tgot_\Qbb$ generated by $\wedge$.

We consider the $K^3$-action on $N:=\T^* K\times \T^* K \times E$. 

\begin{defi}
A non-zero element $(\gamma_1,\gamma_2,\gamma_3)\in\tgot^3$ 
is called {\em admissible} if the elements $\gamma_i$ are rational and if $\dim_{K^3}(N^{(\gamma_1,\gamma_2,\gamma_3)})-\dim_{K^3}(N)\in\{0,1\}$.
\end{defi}

Let $\Rgot$ be the set of roots for $(K,T)$, and let $W=N(T)/T$ be the Weyl group. The set of weights 
for the $T$-action on $E$ is denoted $\Rgot_E$. 
If $\gamma\in\tgot$, we denote by $(\Rgot\cup\Rgot_E)\cap\gamma^\perp$ the subsets of weight vanishing against $\gamma$.

If $w=(w_1,w_2,w_3)\in W^3$ and $\gamma\in\tgot$, we write $\gamma_w=(w_1\gamma,w_2\gamma,w_3\gamma)$. We start with the 
following lemma whose proof is left to the reader.

\begin{lem}
\begin{enumerate}
\item $N^{(\gamma_1,\gamma_2,\gamma_3)}\neq \emptyset$ if and only if $\gamma_1,\gamma_2\in W\gamma_3$.
\item $\dim_{K^3}(N)=\dim_T(\kgot\times E)=\dim(\tgot)-\dim({\rm Vect}(\Rgot\cup\Rgot_E))$.
\item $\dim_{K^3}(N^{\gamma_w})=\dim_T(\kgot^{\gamma}\times E^\gamma)=\dim(\tgot)-\dim({\rm Vect}((\Rgot\cup\Rgot_E)\cap\gamma^\perp))$.
\end{enumerate}
\end{lem}

The following result is a direct consequence of the previous lemma.

\begin{lem}
The admissible elements relative to the $K^3$-action on $\T^* K\times \T^* K \times E$ are of the form $\gamma_w$ where $w\in W^3$ and $\gamma$ is 
a non-zero rational element satisfying ${\rm Vect}(\Rgot\cup\Rgot_E)\cap \gamma^\perp={\rm Vect}((\Rgot\cup\Rgot_E)\cap\gamma^\perp)$.
\end{lem}

\subsection{Ressayre's data}

\begin{defi}
\begin{enumerate}
\item Consider the linear action $\rho: G\to {\rm GL}_\Cbb(V)$ of a compact Lie group on a complex vector space $V$. For any $(\eta,a)\in\ggot\times\Rbb$, we define the vector subspace
$V^{\eta=a}=\{v\in V, d\rho(\eta)v=i av\}$.  Thus, for any $\eta\in\ggot$, we have the decomposition
$V=V^{\eta>0}\oplus V^{\eta=0}\oplus V^{\eta<0}$ where $V^{\eta>0}=\sum_{a>0}V^{\eta=a}$, and $V^{\eta<0}=\sum_{a<0}V^{\eta=a}$.
\item The real number $\tr_{\eta}(V^{\eta>0})$ is defined as the sum $\sum_{a>0}a\,\dim(V^{\eta=a})$.
\end{enumerate}
\end{defi}

We consider an admissible element $\gamma_w=(w_1\gamma,w_2\gamma,w_3\gamma)$. The submanifold fixed  by $\gamma_w$ 
is $N^{\gamma_w}= w_1 K_\Cbb^\gamma w_3^{-1}\times w_2 K_\Cbb^\gamma w_3^{-1}\times E^{w_3\gamma}$. There is a canonical isomorphism of the manifold $N^{\gamma_w}$ equipped with the action of $w_1K^\gamma w_1^{-1}\times w_2K^\gamma w_2^{-1}\times w_3K^\gamma w_3^{-1}$ 
with the manifold $K_\Cbb^\gamma \times K_\Cbb^\gamma \times E^{\gamma}$ equipped with the action of $K^\gamma \times K^\gamma \times K^\gamma$.
The tangent bundle $(\T N\vert_{N^{\gamma_w}})^{\gamma_w>0}$ is isomorphic to 
$N^{\gamma_w}\times \kgot_\Cbb^{\gamma>0}\times \kgot_\Cbb^{\gamma>0}\times E^{\gamma>0}$. 

The choice of positive roots $\Rgot^+$ induces a decomposition $\kgot_\Cbb=\ngot\oplus\tgot_\Cbb\oplus \overline{\ngot}$, where 
$\ngot=\sum_{\alpha\in\Rgot^+}(\kgot\otimes\Cbb)_\alpha$. We consider the map 
$$
\rho^{\gamma,w}: K_\Cbb^\gamma \times K_\Cbb^\gamma\times	 E^\gamma \longrightarrow 
\hom \left(\ngot^{w_1\gamma>0}\times \ngot^{w_2\gamma>0}\times \ngot^{w_3\gamma>0}, \kgot_\Cbb^{\gamma>0}\times \kgot_\Cbb^{\gamma>0}\times E^{\gamma>0}\right)
$$
defined by the relation
$$
\rho^{\gamma,w}(x,y,v): (X,Y,Z)\mapsto ((w_1 x)^{-1}X-w_3^{-1}Z\, ;\, (w_2 y)^{-1}Y-w_3^{-1}Z\, ;\, (w_3^{-1}Z)\cdot v),
$$
for any $(x,y,v)\in K_\Cbb^\gamma \times K_\Cbb^\gamma\times E^\gamma$.
\medskip

\begin{defi}
$(\gamma,w)\in\tgot\times W^3$ is a Ressayre's data if 
\begin{enumerate}
\item $\gamma_w$ is admissible,
\item $\exists (x,y,v)$ such that $\rho^{\gamma,w}(x,y,v)$ is bijective.
\end{enumerate}
\end{defi}

\begin{rem}
In \cite{pep-ressayre-pair}, the Ressayre's data were called {\em regular infinitesimal $B$-Ressayre's pairs}.
\end{rem}

Since the linear map $\rho^{\gamma,w}(x,y,v)$ commutes with the $\gamma$-actions, we obtain the following necessary conditions.

\begin{lem}\label{lem:relations-A-B}
If $(\gamma,w)\in\tgot\times W^3$ is a Ressayre's data, then
\begin{itemize}
\item Relation (A) : $\sum_{i=1}^3\dim (\ngot^{w_i\gamma>0})= 2\dim(\kgot_\Cbb^{\gamma>0})+\dim( E^{\gamma>0})$.
\item Relation (B) : $\sum_{i=1}^3 \tr_{w_i\gamma}(\ngot^{w_i\gamma>0})=2\tr_{\gamma}(\kgot_\Cbb^{\gamma>0})+ \tr_{\gamma}(E^{\gamma>0})$.
\end{itemize}
\end{lem}

\subsection{Cohomological characterization of Ressayre's data}

Let $\gamma\in\tgot$ be a rational element. We denote by $B\subset K_\Cbb$ the Borel subgroup with Lie algebra $\bgot=\tgot_\Cbb\oplus \ngot$. Consider the parabolic subgroup $P_\gamma\subset K_\Cbb$ defined by
\begin{equation}\label{eq:P-gamma}
P_\gamma=\{g\in K_\Cbb, \lim_{t\to\infty}\exp(-it\gamma)g\exp(it\gamma)\ {\rm exists}\}.
\end{equation}
We work with the projective variety $\Fcal_\gamma:=K_\Cbb/P_\gamma$. We associate to any $w\in W$, the Schubert cell
$$
\Xgot^o_{w,\gamma}:= B [w]\subset \Fcal_\gamma,
$$
and the Schubert variety $\Xgot_{w,\gamma}:=\overline{\Xgot^o_{w,\gamma}}$. If $W^\gamma$ denotes the subgroup of $W$ that fixes 
$\gamma$, we see that  the Schubert cell $\Xgot^o_{w,\gamma}$ and the Schubert variety $\Xgot_{w,\gamma}$ depends only of the class of 
$w$ in $W/W^\gamma$.

We consider the cohomology\footnote{Here, we use singular cohomology with integer coefficients.} ring $H^*(\Fcal_\gamma,\Zbb)$ of $\Fcal_\gamma$. 
If $Y$ is an irreducible closed subvariety of $\Fcal_\gamma$, we denote by $[Y]\in H^{2n_Y}(\tilde{\Fcal}_\gamma,\Zbb)$ its cycle class in cohomology : here 
$n_Y={\rm codim}(Y)$. Recall that the cohomology class $[pt]$ associated to a singleton $\{pt\}\subset \Fcal_\gamma$  is a basis of 
$H^{\maxx}(\Fcal_\gamma,\Zbb)$.

To an oriented real vector bundle $\Ecal\to N$ of rank $r$, we can associate its Euler class ${\rm Eul}(\Ecal)\in H^{r}(N,\Zbb)$. 
When $\Ecal\to N$ is a complex vector bundle, then ${\rm Eul}(\Ecal_\Rbb)$ corresponds to the top Chern class $c_{p}(\Ecal)$. 
Here $p$ is the complex rank of $\Ecal$, and $\Ecal_\Rbb$ means $\Ecal$ viewed as a real vector bundle oriented by its complex structure (see \cite{Bott-Tu}, \S 21).

The isomorphism $E^{\gamma>0}\simeq E/E^{\gamma\leq 0}$ shows that $E^{\gamma>0}$ can be viewed as a $P_\gamma$-module. Let 
$\Ecal^{\gamma>0}= K_\Cbb\times_{P_\gamma}E^{\gamma>0}$ be the corresponding complex vector bundle on $\Fcal_\gamma$.
In the following proposition, we denote simply by ${\rm Eul}(E^{\gamma>0})$ the Euler class ${\rm Eul}(\Ecal^{\gamma>0}_\Rbb)\in H^*(\Fcal_\gamma,\Zbb)$.

The following characterization of Ressayre's data was obtained in \cite{pep-ressayre-pair}, \S 6.
\begin{prop}
$(\gamma,w)\in\tgot\times W^3$ is a Ressayre's data if and only if
\begin{enumerate}
\item $\gamma_w$ is admissible,
\item Relation (B) holds,
\item The following relation holds in $H^*(\Fcal_\gamma,\Zbb)$ :
\begin{equation}\label{eq:cohomological}
[\Xgot_{w_1,\gamma}]\cdot [\Xgot_{w_2,\gamma}]\cdot[\Xgot_{w_3,\gamma}]\cdot {\rm Eul}(E^{\gamma>0})=k [pt],\quad k\geq 1.
\end{equation}
\end{enumerate}
\end{prop}

\begin{rem}
Notice that relation (A) is equivalent to $\sum_{i=1}^3{\rm codim}(\Xgot_{w_i,\gamma})+\dim(E^{\gamma>0})=\dim(\Fcal_\gamma)$, hence relation (A) follows from  (\ref{eq:cohomological}).
\end{rem}
\medskip

\subsection{Convex cone $\Delta(K_\Cbb\times K_\Cbb\times E)$ : equations of the facets}

The following result is proved in \cite{pep-ressayre-pair}, \S 6. 

\begin{theo}\label{theo:delta-cas-general} Let $E$ be a $K$-module such that $\Cbb[E]^K=\Cbb$. 
An element $(\xi_1,\xi_2,\xi_3)\in(\tgot^*_{\geq 0})^3$ belongs to $\Delta(K_\Cbb\times K_\Cbb\times E)$ if and only 
\begin{equation}\label{eq:inegalite-admissible}
\langle \xi_1,w_1\gamma\rangle+ \langle \xi_2,w_2\gamma\rangle+\langle \xi_3,w_3\gamma\rangle\geq 0
\end{equation}
for any $(\gamma,w)\in\tgot\times W^3$ that is a Ressayre's data, that is to say satisfying the following properties:
\begin{enumerate}
\item[a)] $\gamma$ is a non-zero rational element.
\item[b)] ${\rm Vect}(\Rgot\cup\Rgot_E)\cap \gamma^\perp={\rm Vect}((\Rgot\cup\Rgot_E)\cap\gamma^\perp)$.
\item[c)] $[\Xgot_{w_1,\gamma}]\cdot [\Xgot_{w_2,\gamma}]\cdot[\Xgot_{w_3,\gamma}]\cdot {\rm Eul}(E^{\gamma>0})=k [pt],\ k\geq 1$ in $H^*(\Fcal_\gamma,\Zbb)$.
\item[d)] Relation (B) holds : $\sum_{i=1}^3 \tr_{w_i\gamma}(\ngot^{w_i\gamma>0})=2\tr_{\gamma}(\kgot_\Cbb^{\gamma>0})+ \tr_{\gamma}(E^{\gamma>0})$.
\end{enumerate}
\end{theo}

\subsection{Remark on the saturation property}

The semigroup $\Delta^\Zbb(K_\Cbb\times K_\Cbb\times E)\subset (\wedge^*_+)^3$ is called {\em saturated} if for any $\theta\in(\wedge^*_+)^3$ and any $N\geq 1$ we have 
$N\theta \in \Delta^\Zbb(K_\Cbb\times K_\Cbb\times E)$ only if $\theta \in \Delta^\Zbb(K_\Cbb\times K_\Cbb\times E)$.

\begin{prop}\label{prop:saturation-remark}
\begin{enumerate}
\item We have 
$$
\Qbb^{>0}\cdot\Delta^\Zbb(K_\Cbb\times K_\Cbb\times E)=\Delta(K_\Cbb\times K_\Cbb\times E)\cap(\tgot^*_\Qbb)^3.
$$
\item The semigroup $\Delta^\Zbb(K_\Cbb\times K_\Cbb\times E)$ is saturated if and only if 
$$
\Delta^\Zbb(K_\Cbb\times K_\Cbb\times E)=\Delta(K_\Cbb\times K_\Cbb\times E)\cap(\wedge^*_+)^3.
$$
\end{enumerate}
\end{prop}

{\em Proof.} Let us prove the first point. The inclusion $\Qbb^{>0}\cdot\Delta^\Zbb(K_\Cbb\times K_\Cbb\times E)\subset\Delta(K_\Cbb\times K_\Cbb\times E)\cap(\tgot^*_\Qbb)^3$ follows from the definitions. Let us explain why the opposite inclusion is a consequence of the $[Q,R]=0$ theorem. 

We consider the proper moment map $\Phi:K_\Cbb\times K_\Cbb\times E\to(\kgot^*)^3$.
For any $\mu=(\mu_1,\mu_2,\mu_3)\in(\wedge^*_+)^3$, we denote by $m_E(\mu)$ the multiplicity of 
$V_\mu=V_{\mu_1}^{K,1}\otimes V_{\mu_2}^{K,2}\otimes V_{\mu_3}^{K,3}$ in the coordinate ring $\Cbb[K_\Cbb\times K_\Cbb\times E]$, 
and we consider the reduced space 
 $$
 \Mcal_\mu:=\Phi^{-1}(K\mu_1\times K\mu_2\times K\mu_2)/K\times K\times K.
 $$
that is equipped with the line bundle
 $$
  \Lcal_{\mu}=\Phi^{-1}(K\mu_1\times K\mu_2\times K\mu_2)\times_{K_{\mu_1}\times K_{\mu_2}\times K_{\mu_2}}
 \left( \Cbb_{-\mu_1}\otimes\Cbb_{-\mu_2}\otimes\Cbb_{-\mu_3}\right).
 $$
 Suppose that $\Mcal_\mu$ is non-empty. Then $\Mcal_{\mu}$ is a complex-projective variety, a projective embedding being 
given by the Kodaira map $\Mcal_\mu \to \Pbb(H^0(\Mcal_\mu, \Lcal_{\mu}^{\otimes k}))$ for all sufficiently large $k$ 
(see Theorem 2.17 in \cite{Sjamaar-95}). 
Moreover, the $[Q,R]$ theorem says that for all $k\geq 1$, we have $m_E(k\mu)=\dim H^0(\Mcal_{\mu}, \Lcal_{\mu}^{\otimes k})$: 
hence $m_E(k\mu)\neq 0$ for $k$ sufficiently large.
 
Let $\xi\in \Delta(K_\Cbb\times K_\Cbb\times E)\cap(\tgot^*_\Qbb)^3$: let $N\geq 1$ such that $\mu_o:=N\xi\in (\wedge^*_+)^3$. 
By definition, the reduced space $\Mcal_{\mu_o}$ is non-empty. So there exists, $k_o\geq 1$ such that $m_E(k_o\mu_o)\neq 0$,  i.e. 
$k_o\mu_o\in \Delta^\Zbb(K_\Cbb\times K_\Cbb\times E)$. 
We have proved that  $\xi=\frac{1}{k_oN}k_o\mu\in \Qbb^{>0}\cdot\Delta^\Zbb(K_\Cbb\times K_\Cbb\times E)$.
 
The first point is settled and the second one is an immediate consequence of the first one. $\Box$

\section{Saturated semigroups}\label{sec:saturated semi-group}

For any $n\geq 1$, we consider the semigroup $\wedge_n^+=\{(\lambda_1\geq \cdots\geq \lambda_n)\}\subset\Zbb^n$ that parametrizes the irreducible representations
 of the unitary group $U(n)$. When $\lambda\in \wedge_n^+$, the notation $\lambda\geq 0$ (resp. $\lambda\leq 0$) means that $\lambda_n\geq 0$ (resp. $\lambda_1\leq 0$), and 
 ${\rm length}(\lambda)$ is the number of non-zero coordinates $\lambda_i$. To $\lambda=(\lambda_1\geq \cdots\geq \lambda_n)\in \wedge_n^+$, we associate $\lambda^*=(-\lambda_n\geq \cdots\geq -\lambda_1)\in \wedge_n^+$ :  the representation $V^{U(n)}_{\lambda^*}$ is then the dual of $V^{U(n)}_{\lambda}$. 
 
Let start by recalling  the properties of the semigroup $\horn^\Zbb(n)$ of $(\wedge_n^+)^3$ defined by the relations
$$
(\lambda,\mu,\nu)\in \horn^{\Zbb}(n)\Longleftrightarrow \left[V^{U(n)}_\nu\,:\,V^{U(n)}_\lambda\otimes V^{U(n)}_\mu\right]\neq 0.
$$
First, the convex cone of $(\Ccal_n)^3$ defined as the closure of $\Qbb^{>0}\cdot \horn^{\Zbb}(n)$ corresponds to  $\horn(n)$ (see Definition 
\ref{defi:horn-cone}). Moreover, thanks to the saturation Theorem of A. Knutson and T. Tao 
\cite{Knutson-Tao-99}, we know that  an element $(\lambda,\mu,\nu)\in (\wedge_n^+)^3$ belongs to the semigroup $\horn^{\Zbb}(n)$ if and only if 
$(\lambda,\mu,\nu)\in \horn(n)$ (see Proposition \ref{prop:saturation-remark}).

In the rest of this section we work with the compact Lie group $K=U(p)\times U(q)$,  so that $K_\Cbb=GL_p\times GL_q$. 
If $\lambda=(\lambda',\lambda'')\in \wedge_p^+\times\wedge_q^+$, then $V_\lambda:=V^{U(p)}_{\lambda'}\otimes V^{U(q)}_{\lambda''}$ denotes the irreducible representation of $U(p)\times U(q)$ with highest weight $\lambda$.  

Recall that for any $p,q\geq 1$,  $M_{p,q}$ denotes the vector space of $p\times q$ complex matrices.

The purpose of this section is the study of the following semigroups of $(\wedge_p^+\times\wedge_q^+)^3$.
\begin{defi} Let $(\lambda,\mu,\nu)\in(\wedge_p^+\times\wedge_q^+)^3$.
\begin{itemize} 
\item The semigroup $\horn^{\Zbb}(p,q)$ is defined by the conditions:
$$
 (\lambda,\mu,\nu)\in \horn^{\Zbb}(p,q)\Longleftrightarrow \left[V_\nu\, :\, V_\lambda\otimes V_\mu\otimes{\rm Sym}(M_{p,q})\right]\neq 0.
$$
\item The semigroup ${\rm Q}^{\Zbb}(p,q)$ is defined by the conditions:
$$
 (\lambda,\mu,\nu)\in {\rm Q}^{\Zbb}(p,q)\Longleftrightarrow  
 \begin{cases}
\lambda\leq 0,\\
\mu\leq 0, \\
 \left[V_\lambda\otimes V_\mu\otimes V_{\nu}\otimes{\rm Sym}(M_{p,q})\right]^{U(p)\times U(q)}\neq 0.
\end{cases}
 $$
\item The semigroup ${\rm S}^{\Zbb}(p,q)$ is defined by the conditions:
$$
 (\lambda,\mu,\nu)\in {\rm S}^{\Zbb}(p,q)\Longleftrightarrow \left[V_\lambda\otimes V_\mu\otimes V_\nu\otimes {\rm Sym}(\Cbb^p\otimes\Cbb^q)\right]^{U(p)\times U(q)}\neq 0.
$$
\end{itemize}
\end{defi}

\subsection{The semigroup ${\rm Q}^\Zbb(p,q)$}

\begin{figure}
\begin{center}
  \includegraphics[width=3 in]{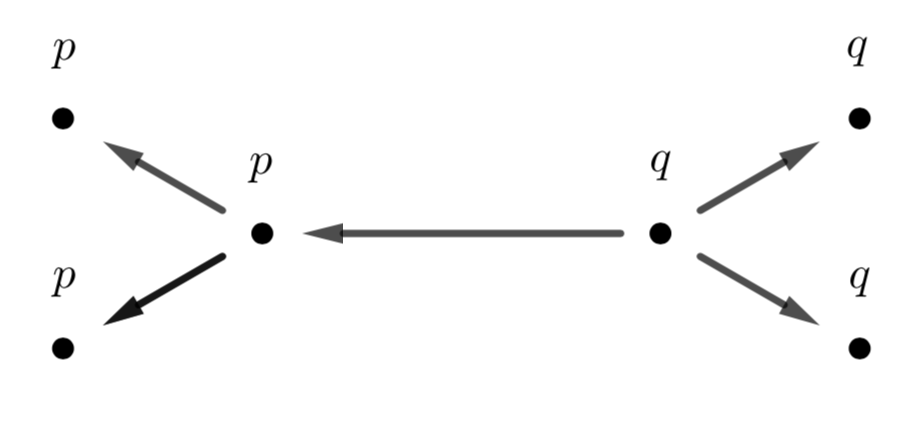}
\caption{Quiver $Q$ with dimension vector $\mathbf{v}_{p,q}$}
 \label{Quiver}
\end{center}
\end{figure}

Let $p,q\geq 1$. We consider the quiver $Q$ of Figure 1, with dimension vector $\mathbf{v}_{p,q}=(p,p,p,q,q,q)$. The vector space 
$$
{\rm Rep}(Q,\mathbf{v}_{p,q})=(M_{p,p}\times M_{q,q})^2\times M_{p,q}
$$
admits a natural action of the algebraic group ${\rm GL}(Q,\mathbf{v}_{p,q})=(GL_p\times GL_q)^3$ that we recall. Take 
$g=(g_1,g_2,g_3)\in {\rm GL}(Q,\mathbf{v}_{p,q})$ with $g_i=(g_i',g_i'')\in GL_p\times GL_q$ and $(X_1,X_2,Y)\in {\rm Rep}(Q,\mathbf{v})$ where 
$X_i=(X'_i,X''_i)\in M_{p,p}\times M_{q,q}$ and $Y\in M_{p,q}$. Then $g\cdot X= (x_1,y_2,y)$  where 
$x_i=(g_i'X'_i(g_3')^{-1},g_i''X''_i(g_3'')^{-1})$ and $y=g_3'Y(g_3'')^{-1}$.

We consider the multipicity map $m: (\wedge^+_p\times\wedge^+_q)^3\to \Nbb$ defined by 
$$
m(\lambda,\mu,\nu)=\dim \left[V_\lambda\otimes V_\mu\otimes V_\nu\otimes{\rm Sym}(M_{p,q}) \right]^{GL_p\times GL_q}.
$$

\begin{lem}
The coordinate ring $\Cbb[{\rm Rep}(Q,\mathbf{v}_{p,q})]$, viewed as ${\rm GL}(Q,\mathbf{v}_{p,q})$-module, admits the following decomposition
$$
\Cbb[{\rm Rep}(Q,\mathbf{v}_{p,q})]=\sum_{\lambda\leq 0,\mu\leq 0, \nu}\,m(\lambda,\mu,\nu)\,V^1_\lambda\otimes V^2_\mu\otimes V^3_\nu.
$$
\end{lem}

{\em Proof.} It is due to the fact that the ${\rm GL}(Q,\mathbf{v}_{p,q})$-module $\Cbb[(M_{p,p}\times M_{q,q})^2]$ admits the decomposition
$\Cbb[(M_{p,p}\times M_{q,q})^2]=\sum_{\lambda\leq 0,\mu\leq 0}V^1_\lambda\otimes V^2_\mu\otimes V^3_{\lambda^*}\otimes V^3_{\mu^*}$.
$\Box$

The previous lemma shows that ${\rm Q}^\Zbb(p,q)$ corresponds to the semigroup of highest weights 
associated to the action of the group ${\rm GL}(Q,\mathbf{v}_{p,q})$ on the coordinate ring  $\Cbb[{\rm Rep}(Q,\mathbf{v}_{p,q})]$. 

\begin{prop}
The semigroup ${\rm Q}^\Zbb(p,q)$ is saturated.
\end{prop}

\begin{figure}
\begin{center}
  \includegraphics[width=5 in]{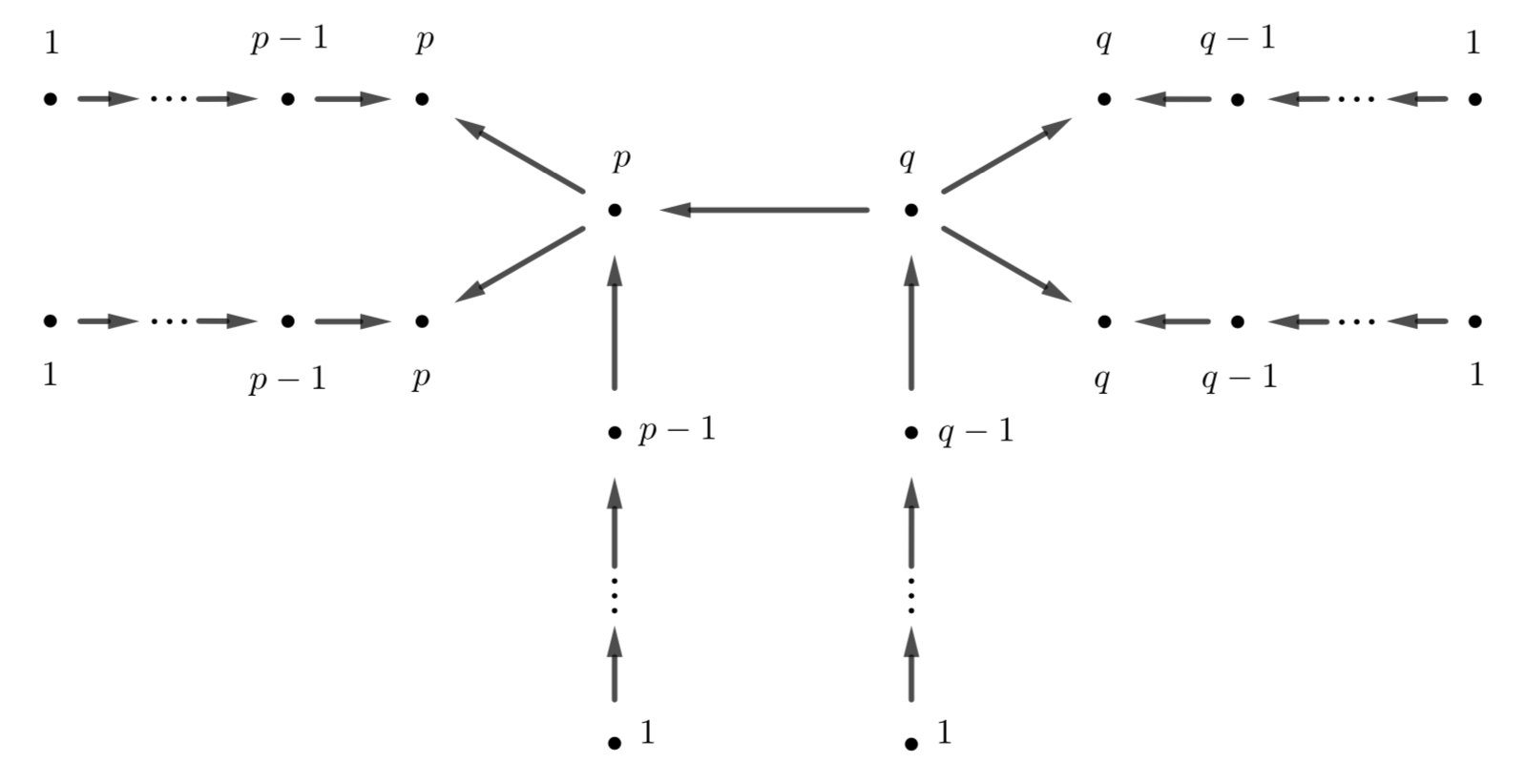}
\caption{Quiver $\widetilde{Q}$ with dimension vector $\widetilde{\mathbf{v}}_{p,q}$}
 \label{Quiver-tilde}
\end{center}
\end{figure}

{\em Proof.} This is a direct consequence of the Derksen-Weyman saturation theorem \cite{Derksen-Weyman}, which asserts that, for a quiver without cycles, 
the semigroup of weights of semi-invariants is saturated. Indeed, augment the quiver $(Q,\mathbf{v}_{p,q})$ to the quiver $(\widetilde{Q},\widetilde{\mathbf{v}}_{p,q})$ 
(see Figure \ref{Quiver-tilde}). Then, using the Cauchy formula for the decomposition of $\otimes_{k=1}^{n-1}\Cbb[(\Cbb^k)^*\otimes\Cbb^{k+1}]$ under the action 
of $\prod_{k=1}^n GL_k$, one sees that there is a bijective morphism between the semigroup of weights of semi-invariants of the coordinate ring 
$\Cbb[{\rm Rep}(\widetilde{Q},\widetilde{\mathbf{v}}_{p,q})]$ under the action of $(\prod_{k=1}^p GL_k)^3\times(\prod_{\ell=1}^q GL_\ell)^3$ and the semigroup 
${\rm Q}^\Zbb(p,q)$. $\Box$

\subsection{The semigroups $\horn^\Zbb(p,q)$ and ${\rm S}^{\Zbb}(p,q)$}

We use the involution $\Theta: (\wedge_p^+\times\wedge_q^+)^3\to (\wedge_p^+\times\wedge_q^+)^3$ that sends $(\lambda,\mu,\nu)$ to 
$$
\left(\begin{pmatrix}
\lambda' \\
(\lambda'')^*
\end{pmatrix},
\begin{pmatrix}
\mu' \\
(\mu'')^*
\end{pmatrix},
\begin{pmatrix}
(\nu')^* \\
\nu''
\end{pmatrix}
\right).
$$

Let us denote by ${\bf 1}\in \wedge_p^+\times\wedge_q^+$ the vector $(1,\ldots,1,1,\ldots,1)$. 
The main purpose of this section is the following result.

\begin{prop}\label{prop:sturation-S-Horn}
\begin{enumerate}
\item For any $(\lambda,\mu,\nu)\in \horn^\Zbb(p,q)$, there exists $k_o\geq 0$ such that 
$(\lambda-k\mathbf{1},\mu-k\mathbf{1},\nu^*+2k\mathbf{1})\in {\rm Q}^\Zbb(p,q)$, $\forall k\geq k_o$.
\item $(\lambda,\mu,\nu)\in \horn^\Zbb(p,q)$ if and only if $\Theta(\lambda,\mu,\nu)\in {\rm S}^{\Zbb}(p,q)$.
\item The semigroups $\horn^\Zbb(p,q)$ and ${\rm S}^{\Zbb}(p,q)$ are saturated.
\end{enumerate}
\end{prop}

{\em Proof.} Let $(\lambda,\mu,\nu)\in \horn^\Zbb(p,q)$ :  it means that \break
$\left[V_\lambda\otimes V_\mu\otimes V_{\nu^*}\otimes{\rm Sym}(M_{p,q}) \right]^{GL_p\times GL_q}\neq 0$ and so 
$(\lambda,\mu,\nu^*)\in {\rm Q}^\Zbb(p,q)$ if $\lambda\leq 0$ and $\mu\leq 0$. We notice that 
$(\lambda+k{\bf 1},\mu+k{\bf 1},\nu+2k{\bf 1})\in \horn^\Zbb(p,q)$, $\forall k\in \Zbb$. Let $k_o=\sup(|\lambda_1|,|\mu_1|)$: we see that 
$\lambda-k{\bf 1}\leq 0$ and $\mu-k{\bf 1}\leq 0$ if $k\geq k_o$, and consequently 
$(\lambda-k{\bf 1},\mu-k{\bf 1},\nu^*+2k{\bf 1})\in {\rm Q}^\Zbb(p,q)$ if $k\geq k_o$. The first point is settled.

Cauchy formulas give the decompositions 
\begin{align*}
{\rm Sym}(M_{p,q})={\rm Sym}( \Cbb^p\otimes(\Cbb^q)^*)&= 
\sum_{\stackrel{{\rm length}(a)\leq \inf(p,q)}{a\geq 0}}V^{GL_p}_{a}\otimes V^{GL_q}_{a^*},\\
{\rm Sym}( \Cbb^p\otimes\Cbb^q)&= 
\sum_{\stackrel{{\rm length}(a)\leq \inf(p,q)}{a\geq 0}}V^{GL_p}_{a}\otimes V^{GL_q}_{a}.
\end{align*}

For the second point, we have to compare the following cases :
\begin{itemize}
\item $(\lambda,\mu,\nu)\in \horn^\Zbb(p,q)$ if there exists $a\geq 0$ with ${\rm length}(a)\leq \inf(p,q)$ such that 
both conditions hold
\begin{align*}
\left[V_{\lambda'}\otimes V_{\mu'}\otimes V_{(\nu')^*}\otimes V_a \right]^{GL_p}     &\neq  0\\
\left[V_{\lambda''}\otimes V_{\mu''}\otimes V_{(\nu'')^*}\otimes V_{a^*}\right]^{GL_q} &\neq  0.
\end{align*}
\item $(\lambda,\mu,\nu)\in {\rm S}^\Zbb(p,q)$ if there exists $a\geq 0$ with ${\rm length}(a)\leq \inf(p,q)$ such that 
both conditions hold
\begin{align*}
\left[V_{\lambda'}\otimes V_{\mu'}\otimes V_{\nu'}\otimes V_a \right]^{GL_p}     &\neq  0\\
\left[V_{\lambda''}\otimes V_{\mu''}\otimes V_{\nu''}\otimes V_{a}\right]^{GL_q} &\neq  0.
\end{align*}
\end{itemize}
It is then immediate to conclude that $(\lambda,\mu,\nu)\in \horn^\Zbb(p,q)$ if and only if  
$\Theta(\lambda,\mu,\nu)\in {\rm S}^\Zbb(p,q)$.

Let us check that $\horn^\Zbb(p,q)$ is saturated. Let $(\lambda,\mu,\nu)$ such that \break $N(\lambda,\mu,\nu)\in \horn^\Zbb(p,q)$ for some $N\geq 1$. 
If we use the first point, we know that there exists $k_o\geq 0$, such that 
$$
N(\lambda,\mu,\nu^*)- k({\bf 1},{\bf 1},-2\cdot{\bf 1})\in {\rm Q}^\Zbb(p,q),\quad \forall k\geq k_o.
$$
Take $k=Nk_o\geq k_o$ : we obtain that $N(\lambda-k_o{\bf 1},\mu-k_o{\bf 1},\nu^*+2k_o{\bf 1})\in {\rm Q}^\Zbb(p,q)$. It follows that 
$(\lambda-k_o{\bf 1},\mu-k_o{\bf 1},\nu^*+2k_o{\bf 1})\in {\rm Q}^\Zbb(p,q)$ because the semigroup  ${\rm Q}^\Zbb(p,q)$ is saturated. 
Finally, the last relation implies that $(\lambda,\mu,\nu)\in \horn^\Zbb(p,q)$. 

We have verified that the semigroup $\horn^\Zbb(p,q)$ is saturated, and the second point also allows us to conclude 
that the semigroup ${\rm S}^\Zbb(p,q)$ is saturated. $\Box$

\subsection{Final remarks}

We have seen that the semigroups $\horn^\Zbb(p,q)$, ${\rm S}^\Zbb(p,q)$ and ${\rm Q}^\Zbb(p,q)$ are all related. 
Thus, the associated convex cones $\horn(p,q)=\overline{\Qbb^{>0}\cdot\horn^\Zbb(p,q)}$, ${\rm S}(p,q)=\overline{\Qbb^{>0}\cdot{\rm S}^\Zbb(p,q)}$ and 
${\rm Q}(p,q)=\overline{\Qbb^{>0}\cdot{\rm Q}^\Zbb(p,q)}$ are also interdependent. The calculation of one entails those of the others. 
In this paper, we obtain a recursive description of the Horn cone $\horn(p,q)$ through the calculation of ${\rm S}(p,q)$. 

Let $Q_o$ be a quiver without cycle which is equipped with a dimension vector ${\bf v_o}$. In \cite{Baldoni-Vergne-Walter19}, 
V. Baldoni, M. Vergne and M. Walter have proposed a recursive description of the cone generated by the highest weights 
associated to the action of the group ${\rm GL}(Q_o,\mathbf{v_o})$ on the coordinate ring  $\Cbb[{\rm Rep}(Q_o,\mathbf{v_o})]$. 
By applying their result to the quiver $(Q,\mathbf{v}_{p,q})$ in Figure \ref{Quiver}, this should also allow a 
recursive description of the Horn cone $\horn(p,q)$.

\section{Convex cone ${\rm S}(p,q)$}\label{sec:polyhedron-s-p-q}

In this section, we apply the results of \S \ref{sec:K-trois-E} to the case where $K_\Cbb=GL_p\times GL_q$ and the $K_\Cbb$-module is 
$E=\Cbb^p\otimes\Cbb^q$. The  coordinate ring $\Cbb[(GL_p\times GL_q)^2\times \Cbb^p\otimes\Cbb^q]$, viewed as a $(GL_p\times GL_q)^3$-module, admits the following decomposition 
$$
\Cbb[(GL_p\times GL_q)^2\times \Cbb^p\otimes\Cbb^q]=\sum_{\lambda,\mu,\nu\in\wedge_p^+\times\wedge_q^+}m(\lambda,\mu,\nu) \ V_\lambda^1 \otimes V_\mu^2 \otimes V_\nu^3,
$$
where $m(\lambda,\mu,\nu)  =\dim [V_\lambda\otimes V_\mu\otimes V_\nu\otimes {\rm Sym}(\Cbb^p\otimes\Cbb^q)]^{GL_p\times GL_q}$.

We see that the semigroup $\Delta^\Zbb((GL_p\times GL_q)^2\times\Cbb^p\otimes\Cbb^q)$ corresponds to ${\rm S}^\Zbb(p,q)$. 
Hence we will denote by ${\rm S}(p,q)$ the convex cone $\Delta((GL_p\times GL_q)^2\times\Cbb^p\otimes\Cbb^q)$ : it corresponds to 
the Kirwan polyhedron relative to the Hamiltonian action of $U(p)\times U(q)$ on $(GL_p\times GL_q)^2 \times \Cbb^p\otimes \Cbb^q$.

%

\subsection{Admissible elements for $(GL_p\times GL_q)^2 \times \Cbb^p\otimes \Cbb^q$}

Let $T\simeq U(1)^p\times U(1)^q$ be the maximal torus of $K=U(p)\times U(q)$ formed by the diagonal matrices. The Lie algebra $\tgot$ admits a canonical identification 
with $\Rbb^p\times\Rbb^q$ through the map ${\bf d}$ and the Weyl group of $(K,T)$ is isomorphic to $W=\Sgot_p\times\Sgot_q$.

The set of roots for $(K,T)$ is $\Rgot=\{\epsilon_i-\epsilon_j, 1\leq i\neq j\leq p\}\cup$\break  $\{\epsilon'_k-\epsilon'_l, 1\leq k\neq l\leq q\}$. The set of weights 
for the $T$-action on $E=\Cbb^p\otimes\Cbb^q$ is $\Rgot_E=\{\epsilon_i+\epsilon'_k, 1\leq i\leq p, 1\leq k\leq q\}$. Let us denote 
$\Rgot_o=\Rgot\cup\Rgot_E$. We first notice that $\Rgot_o^\perp=\Rbb \gamma^0_0$, with $\gamma^0_0={\bf 1}_p\oplus -{\bf 1}_q\in\tgot$.

\begin{defi} For any $(r,s)\in \{0,\ldots,p\}\times  \{0,\ldots,q\}$, we define $\gamma^{r}_{s}=\gamma^{r}\oplus \gamma_{s}\in\tgot$
where $\gamma^{r}=(\underbrace{0,\ldots,0}_{r\ {\rm times}},1,\ldots,1)$ and 
$\gamma_{s}=(-1,\ldots,-1,\underbrace{0,\ldots,0}_{s\ {\rm  times}})$.
\end{defi}

\begin{lem}
Let $\gamma\in\tgot$ be a non-zero rational element such that ${\rm Vect}(\Rgot_o)\cap \gamma^\perp={\rm Vect}(\Rgot_o\cap \gamma^\perp)$. There exists 
 $(r,s)\notin\{ (0,0)\, , \, (p,q)\}$, $w\in W$ and $(a,b)\in\Qbb^{\geq 0}\times\Qbb$ such that $\gamma=a(w\gamma^r_s)+b\gamma^0_0$.
\end{lem}
{\em Proof.} For any $t\in \Rbb$, we define 
$$
\gamma(t)=\sum_{\stackrel{1\leq i\leq p}{\gamma_i=t}} e_i-\sum_{\stackrel{1\leq k\leq q}{\gamma_k=-t}}e'_k.
$$
We notice that $\gamma(t)$ is orthogonal to $\Rgot_o\cap \gamma^\perp$. If ${\rm Vect}(\Rgot_o)\cap \gamma^\perp=$ \break 
${\rm Vect}(\Rgot_o\cap \gamma^\perp)$ holds we get  that $\gamma(t)\in ({\rm Vect}(\Rgot_o)\cap \gamma^\perp)^\perp=\Rbb \gamma^0_0 +\Rbb \gamma$ for all $t\in\Rbb$. Take $t_o$ such that $\gamma(t_o)\neq 0$. Two situations holds.
\begin{enumerate}
\item $\gamma(t_o)\in\Rbb\gamma^0_0$. This case only occurs if $\gamma\in\Qbb\gamma^0_0$.
\item $\gamma(t_o)\notin\Rbb\gamma^0_0$. Then there exists $(r,s)\notin\{ (0,0)\, , \, (p,q)\}$, $w\in W$ and 
$(x,y)\in\Qbb-\{0\}\times\Qbb$ such that $\gamma(t_o)=w\gamma^r_s$ and $\gamma(t_o)=x\gamma+y\gamma^0_0$: hence
$\gamma=\tfrac{1}{x}(w\gamma^r_s)-\tfrac{y}{x}\gamma^0_0$. If $\tfrac{1}{x}>0$, the proof is completed. If $\tfrac{1}{x}<0$, 
we use that $-\gamma^r_s= w_o\gamma^{p-r}_{q-s}-\gamma^0_0$ for some $w_o\in W$ in order to come back to the previous case. $\Box$
\end{enumerate}

In order to describe the facets of the convex cone ${\rm S}(p,q)$, we must consider the following admissible elements: 
\begin{itemize}
\item $\pm(\gamma^0_0,\gamma^0_0,\gamma^0_0)$,
\item $(w_1\gamma^r_s,w_2\gamma^r_s,w_3\gamma^r_s)$ where $w_1,w_2,w_3\in W$ and $(r,s)\notin\{ (0,0)\, , \, (p,q)\}$.
\end{itemize}

\subsection{Admissible elements $\pm(\gamma^0_0,\gamma^0_0,\gamma^0_0)$}

The admissible elements $\pm(\gamma^0_0,\gamma^0_0,\gamma^0_0)$ act on $(GL_p\times GL_q)^2 \times \Cbb^p\otimes \Cbb^q$ trivially. Hence 
$\pm(\gamma^0_0,\gamma^0_0,\gamma^0_0)$ are Ressayre's data, and inequalities  (\ref{eq:inegalite-admissible})  are 
$\pm\left(\langle A,\gamma^0_0\rangle+\langle B,\gamma^0_0\rangle+\langle C,\gamma^0_0\rangle\right)\geq 0$. In other words,  
\begin{equation}\label{eq:gamma-0-0}
\boxed{|\,A'\,| + |\,B'\,| + |\,C'\,| = |\,A''\,| + |\,B''\,| + |\,C''\,|}.
\end{equation}

\subsection{Admissible element $(w_1\gamma^r_s,w_2\gamma^r_s,w_3\gamma^r_s)$}

Recall the relations that a Ressayre's data $(\gamma,w_1,w_2,w_3)$ must satisfy (see Lemma \ref{lem:relations-A-B}):
\begin{align*}
{\rm Relation}\  (A) & : &\sum_{i=1}^3\dim (\ngot^{w_i\gamma>0})   = &\  2\dim(\kgot_\Cbb^{\gamma>0})+\dim( (\Cbb^p\otimes\Cbb^q)^{\gamma>0}),\\
{\rm Relation}\ (B)  & : &\sum_{i=1}^3 \tr_{w_i\gamma}(\ngot^{w_i\gamma>0})  = & \ 2\tr_{\gamma}(\kgot_\Cbb^{\gamma>0})+ \tr_{\gamma}((\Cbb^p\otimes\Cbb^q)^{\gamma>0}).
\end{align*}

It is immediate to see that for $(\gamma^r_s,w_1,w_2,w_3)$, Relations (A) and (B) are equivalent.

We associate to $(w_1,w_2,w_3)\in (\Sgot_p\times\Sgot_q)^3$ the following  subsets :
\begin{itemize}
\item Those of cardinal $r$ : $I'\!=\!w'_1(\{1,\ldots,r\})$, $J'\!=\!w'_2(\{1,\ldots,r\})$, and $K'\!=\!w'_3(\{1,\ldots,r\})$.
\item Those of cardinal $s$ : $I''\!=\!w''_1(\{q-\!s+1,\ldots,q\})$, $J''\!=\!w''_2(\{q-\!s+1,\ldots,q\})$, and $K''\!=\!w''_3(\{q-\!s+1,\ldots,q\})$. 
\end{itemize}

Inequality (\ref{eq:inegalite-admissible}) becomes
$|A'|_{(I')^c} + |B'|_{(J')^c} + |C'|_{(K')^c}\geq |A''|_{(I'')^c} + |B''|_{(J'')^c} + |C''|_{(K'')^c}$ which is equivalent to
\begin{equation}\label{eq:equation-facet}
\boxed{|\,A'\,|_{I'} + |\,B'\,|_{J'} + |\,C'\,|_{K'}\leq |\,A''\,|_{I''} + |\,B''\,|_{J''} + |\,C''\,|_{K''}},
\end{equation}
thanks to (\ref{eq:gamma-0-0}).

\subsection{Schubert classes}

For any $m,n\geq 0$, let $\Gbb(m,n)$ denote the Grassmannian of complex $m$-dimensional linear subspaces of $\Cbb^{m+n}$. The singular cohomology of 
$\Gbb(m,n)$ with integers coefficients is denoted $H^*(\Gbb(m,n),\Zbb)$.

Let $m,n\geq 1$. When a partition $\lambda$ is included in a $m\times n$ rectangle, we write $\lambda\subset m\times n$ : 
$n\geq \lambda_1\geq\cdots\geq \lambda_m\geq 0$.

Denote the set of cardinality $m$ subsets $I=\{i_1<i_2<\cdots<i_m\}$ of $[m+n]=\{1,\ldots,m+n\}$ by $\Pcal^{m+n}_m$. To each 
$I\in \Pcal^{m+n}_m$ we associate $\lambda(I)=(\lambda_1\geq\lambda_2\cdots\geq \lambda_m) \subset m\times n$  
where $\lambda_a= n+a-i_a$ for $a\in [m]$. The map $I \mapsto \lambda(I)$ 
is one to one map between $\Pcal^{m+n}_m$ and the set of partitions of size $m\times n$. The inverse map is denoted by $\lambda\subset m\times n\mapsto I(\lambda)\in \Pcal^{m+n}_m$.

\medskip

We work with the flag $0\subset \Cbb\subset \Cbb^2\subset\cdots\subset \Cbb^{n+m-1}\subset \Cbb^{n+m}$. 
For any partition $\lambda\subset m\times n$, we define the Schubert cell
$$
\Xgot_\lambda^o=\{F\in \Gbb_{m,n},\ \Cbb^{k-1}\cap F\neq \Cbb^{k}\cap F\ \ {\rm if \ and\ only\ if}\ \ k\in I(\lambda)\}.
$$
and the Schubert variety $\Xgot_\lambda=\overline{\Xgot_\lambda^o}$. When the partition is $(k,0,\ldots,0)$ with $1\leq k\leq n$, the corresponding 
Schubert variety is denoted by $\Xgot_{k}$. 
 
Let $B_{n+m}\subset GL_{n+m}$  be the Borel subgroup formed by the upper-triangular matrices. The following facts are well-known (see e.g. 
\cite{Griffiths-Harris,Fulton-97,Manivel-98}) :
\begin{enumerate}
\item $\Xgot_\lambda^o=B_{m+n}\cdot {\rm Vect}(e_i,i\in I(\lambda))$,
\item $\Gbb_{m,n}=\bigcup_{\lambda\subset m\times n}\Xgot_\lambda^o$,
\item ${\rm codim}(\Xgot_\lambda^o)=|\lambda|=\sum_k  \lambda_k$.
\end{enumerate}

Since the Schubert cells define a complex cellular decomposition of the Grassmannian, an immediate consequence is that the fundamental class 
of the Schubert varieties, the {\em Schubert classes}  $\sigma_\lambda=[\Xgot_\lambda]\in H^{2|\lambda|}(\Gbb_{m,n},\Zbb)$, where $\lambda\subset m\times n$, 
form a basis of the cohomology with integers coefficients :
$$
H^*(\Gbb(m,n),\Zbb)=\bigoplus_{\lambda\subset m\times n}\, \Zbb\sigma_\lambda.
$$

When $\lambda=(k,0,\ldots,0)$, we denote by $\sigma_k\in H^{2k}(\Gbb_{m,n},\Zbb)$ the corresponding Schubert class.

We finish this section by recalling that the Grassmannian $\Gbb_{m,n}$ admits the following complex vector bundles :
\begin{enumerate}
\item A canonical vector bundle of rank $m$ : $\Ebb_{m,n}$.
\item A vector bundle of rank $n$, denoted $\Ebb^\perp_{m,n}$,  such that 
$\Ebb_{m,n}\oplus\Ebb^\perp_{m,n}$ is a trivial bundle of rank $m+n$.
\end{enumerate}

\subsection{Cohomological conditions}

Theorem \ref{theo:delta-cas-general} tells us that an element $(A,B,C)\in(\Ccal_p\times\Ccal_q)^3$ belongs to ${\rm S}(p,q)=\Delta((GL_p\times GL_q)^2\times \Cbb^p\otimes\Cbb^q)$ if and only (\ref{eq:gamma-0-0}) holds and (\ref{eq:equation-facet}) holds for any $(w_1,w_2,w_3)\in (\Sgot_p\times\Sgot_q)^3$ and 
any couple $(r,s)\in \{0,\ldots,p\}\times \{0,\ldots,q\}-\{(p,q), (0,0)\}$, satisfying the relation
\begin{equation}\label{eq:cohomological-condition-r-s}
[\Xgot_{w_1,\gamma^r_s}]\cdot [\Xgot_{w_2,\gamma^r_s}]\cdot[\Xgot_{w_3,\gamma^r_s}]\cdot {\rm Eul}(\Vcal^r_s)=k [pt],\ k\geq 1
\end{equation}
in $H^*(\Fcal_{\gamma^r_s},\Zbb)$. 

\medskip

Let us detailed (\ref{eq:cohomological-condition-r-s}). We fix $(r,s)\notin\{ (0,0)\, , \, (p,q)\}$. The parabolic subgroup $P_{\gamma^r_s}\subset GL_p\times GL_q$ associated to $\gamma^r_s$ by (\ref{eq:P-gamma}) 
 is equal to $P_{\gamma^r}\times P_{\gamma_s}$ where 
 $$
 P_{\gamma^r}=
\left(\begin{array}{@{}c|c@{}}
  GL_{r}& *\\ \hline
  0& GL_{p-r}
  \end{array}\right)\subset GL_p
  \quad {\rm and}\quad 
   P_{\gamma_s}=
\left(\begin{array}{@{}c|c@{}}
  GL_{q-s}& *\\ \hline
  0& GL_{s}
  \end{array}\right)\subset GL_q .
$$

Let $\Cbb^{r}\subset \Cbb^p$ and $\Cbb^{q-s}\subset\Cbb^q$ denote respectively the subspaces \break ${\rm Vect}(e_i, 1\leq i\leq r)$ and 
${\rm Vect}(e_j, 1\leq j\leq q-s)$. We use on $\Cbb^p$ and $\Cbb^q$ the canonical bilinear forms $(x,y)\mapsto \sum_{k}x_k y_k$.
Then $(\Cbb^{r})^\perp\subset \Cbb^p$ and $(\Cbb^{q-s})^\perp\subset\Cbb^q$ are respectively the subspaces ${\rm Vect}(e_i, r+1\leq i\leq p)$ and 
${\rm Vect}(e_j, q-s+1\leq j\leq q)$.

The flag variety $\Fcal_{\gamma^r_s}=GL_p/P_{\gamma^r}\times GL_q/P_{\gamma_s}$ admits a canonical identification with 
$\Gbb(r,p-r)\times \Gbb(q-s,s)$ through the map 
$$
([g],[h])\in GL_p/P_{\gamma^r}\times GL_q/P_{\gamma_s}\longmapsto \left(g(\Cbb^{r}),h(\Cbb^{q-s})\right)\in \Gbb(r,p-r)\times \Gbb(q-s,s).
$$

Let $B_p\subset GL_p$ and $B_q\subset GL_q$ be the Borel subgroups formed by the upper-triangular matrices. For any 
$w=(w',w'')\in \Sgot_p\times\Sgot_q$, we consider the Schubert cell
$\Xgot^o_{w,\gamma^r_s}=B_p[w']\times B_q[w'']\subset \Gbb(p-r,r)\times \Gbb(s,q-s)$
and the Schubert variety 
$$
\Xgot_{w,\gamma^r_s}=\overline{B_p[w']}\times \overline{B_q[w'']}=\Xgot_{\mu'}\times \Xgot_{\mu''}
$$
where $\mu'=\lambda(w'\{1,\ldots,r\})\subset r\times p-r$ and $\mu''=\lambda(w''\{1,\ldots,q-s\})\subset q-s\times s$.

If we associate to $(w_1,w_2,w_3)\in (\Sgot_p\times\Sgot_q)^3$, the subsets 
\begin{itemize}
\item $I'\!=\!w'_1(\{1,\ldots,r\})$, $J'\!=\!w'_2(\{1,\ldots,r\})$, and $K'\!=\!w'_3(\{1,\ldots,r\})$,
\item $I''\!=\!w''_1(\{q-\!s+1,\ldots,q\})$, $J''\!=\!w''_2(\{q-\!s+1,\ldots,q\})$, and $K''\!=\!w''_3(\{q-\!s+1,\ldots,q\})$, 
\end{itemize}
the term $[\Xgot_{w_1,\gamma^r_s}]\cdot [\Xgot_{w_2,\gamma^r_s}]\cdot[\Xgot_{w_3,\gamma^r_s}]\in H^*(\Gbb(r,p-r)\times \Gbb(q-s,s),\Zbb)$ 
is then equal to the tensor product of the cohomology classes\footnote{$X^c$ denotes the complement of $X$.}
\begin{align*}
\sigma_{\lambda(I')}\cdot \sigma_{\lambda(J')}\cdot \sigma_{\lambda(K')} & \in H^*(\Gbb(r,p-r),\Zbb),\\
\sigma_{\lambda((I'')^c)}\cdot \sigma_{\lambda((J'')^c)}\cdot \sigma_{\lambda((K'')^c)} &\in H^*(\Gbb(q-s,s),\Zbb).
\end{align*}

The subspace $(\Cbb^p\otimes\Cbb^q)^{\gamma^r_s>0}$ is equal to $(\Cbb^{r})^{\perp}\otimes(\Cbb^{q-s})^\perp$. Hence the vector bundle $\Vcal^r_s$ is 
equal to the tensor product $\Ebb_{r,p-r}^\perp\boxtimes \Ebb_{q-s,s}^\perp$. Let ${\rm Eul}(\Vcal^r_s)\in H^{2(p-r)s}\left(\Gbb(r,p-r)\times \Gbb(q-s,s),\Zbb\right)$ be its Euler class. 

Finally, the cohomological condition (\ref{eq:cohomological-condition-r-s}) says that the product 
$$
\left(\sigma_{\lambda(I')}\otimes \sigma_{\lambda((I'')^c)}\right)\cdot \left(\sigma_{\lambda(J')}\otimes \sigma_{\lambda((J'')^c)}\right)\cdot 
\left(\sigma_{\lambda(K')}\otimes \sigma_{\lambda((K'')^c)}\right)\cdot {\rm Eul}(\Vcal^r_s)
$$
 is a non zero multiple of  $[pt]\in H^{\max}\left(\Gbb(r,p\!-\! r)\!\times\! \Gbb(q\!-\! s,s),\Zbb\right)$.

\section{Computation of the Euler class ${\rm Eul}(\Vcal^r_s)$}\label{sec:euler-class}

Before giving a formula for the Euler class ${\rm Eul}(\Vcal^r_s)$, we need to recall some well-known facts.

\subsection{Polynomial representations}\label{sec:polynomial-rep}

We are concerned here with the polynomial representations of $GL_m$. When a representation $\pi_V: GL_m\to GL(V)$ is polynomial, its character
is an invariant polynomial $\chi_V\in\Cbb[M_{m,m}]$. We denote then by $s_V$ the restriction of $\chi_V$ to the diagonal matrices.

Let $R_+(GL_m)$ denotes the polynomial representation ring of $GL_m$, and let $\bigwedge_m=\Zbb[x_1,\ldots,x_m]^{\Sgot_m}$ be the ring of 
symmetric polynomials, with integral coefficients, in $m$ variables. The map $V\in R_+(GL_m)\mapsto s_V\in\bigwedge_m$ is a ring isomorphism.

For any partition $\lambda$ of length $m$, we associate the irreducible polynomial representation\footnote{When the group is understood, we use the notation $V_\lambda$.} 
$V_{\lambda}^{GL_m}$ of the group $GL_m$ and the Schur polynomial ${\bf s}_\lambda:=s_{V_\lambda}\in\bigwedge_m$. Recall that the Schur polynomials 
${\bf s}_\lambda$ determine a $\Zbb$-basis of $\bigwedge_m$.

We recall the following classical fact (for a proof see \S 3.2.2 in \cite{Manivel-98}).
\begin{theo}\label{theo:phi-n-m}
The map $\phi_{m,n}:{\bigwedge}_m\longrightarrow H^{*}(\Gbb_{m,n},\Zbb)$ defined by the relations 
$$
\phi_{m,n}({\bf s}_\lambda)=
\begin{cases}
\sigma_\lambda\hspace{8mm} {\rm if}\quad \lambda_1\leq n,\\
0\hspace{1cm}{\rm if}\quad \lambda_1> n.
\end{cases}
$$
is a ring morphism. 
\end{theo}

\begin{rem}
Since $V\in R_+(GL_m)\mapsto s_V\in\bigwedge_m$ is a ring isomorphism, we will also denote by 
$\phi_{m,n}:R_+(GL_m)\longrightarrow H^{*}(\Gbb_{m,n},\Zbb)$ the ring morphism 
$V\mapsto \phi_{m,n}(s_V)$.
\end{rem}

If $k\geq 1$, we denote by $1^k$ the partition $(1,\ldots,1,0,\ldots,0)$ where there are $k$-times $1$.

\begin{exam}
\begin{align*}
\phi_{m,n}({\rm Sym}^k(\Cbb^m)) & =
\begin{cases}
\sigma_k\hspace{8mm} {\rm if}\quad 1\leq k \leq n,\\
0\hspace{1cm}{\rm if}\quad k>n.
\end{cases}
\\
\phi_{m,n}({\bigwedge}^k\Cbb^m) & =
\begin{cases}
\sigma_{1^k}\hspace{6mm} {\rm if}\quad 1\leq k \leq m,\\
0\hspace{1cm}{\rm if}\quad k>m.
\end{cases}
\end{align*}
\end{exam}

\subsection{Duality I }\label{sec:duality-1}

We associate to a partition $\lambda\subset m\times n$ it's complementary partition $\widehat{\lambda}\subset m\times n $~:
$\widehat{\lambda}_k=n-\lambda_{m+1-k},\ 1\leq k\leq m$. Recall that 
\begin{itemize}
\item  If $\lambda=\lambda(I)$ then $\widehat{\lambda}=\lambda(\tI)$ where $\tI=\{n+m+1-i\, ;\, i\in I\}$.
\item  If $V_\lambda$ is the irreducible polynomial representation of $GL_m$ associated to $\lambda\subset m\times n$, then 
$(V_\lambda)^*=V_{\widehat{\lambda}}\otimes {\det}^{-n}$.
\end{itemize}

The cohomology class $[pt]\in H^{2nm}(\Gbb_{m,n},\Zbb)$ of top degree associated to a singleton $\{pt\}$ is a basis of 
$H^{2nm}(\Gbb_{m,n},\Zbb)$.

We recall the following classical fact (for a proof see \S 3.2.2 in \cite{Manivel-98}).
\begin{prop}\label{prop:lambda-hat-dual}
Let $\lambda',\lambda\subset m\times n$ be two partitions such that $|\lambda|+|\lambda'|=nm$. Then, the following relations hold in $H^*(\Gbb_{m,n},\Zbb)$ :
$$
\sigma_\lambda\cdot\sigma_{\lambda'}=
\begin{cases}
[pt]\quad {\rm if} \quad \lambda'=\widehat{\lambda},\\
0\hspace{8mm} {\rm if} \quad \lambda'\neq\widehat{\lambda}.
\end{cases}
$$
\end{prop}

The next corollary follows from Theorem \ref{theo:phi-n-m} and Proposition \ref{prop:lambda-hat-dual}.

\begin{coro} Let $\lambda_1,\lambda_2,\lambda_3\subset m\times n$. The following assertions are equivalent :
\begin{itemize}
\item $\sigma_{\lambda_1}\cdot \sigma_{\lambda_2}\cdot\sigma_{\lambda_3}= k[pt],\ k\geq 1$ in $H^{*}(\Gbb_{m,n},\Zbb)$.
\item $\left[V_{\widehat{\lambda}_3} : V_{\lambda_1}\otimes V_{\lambda_2}\right]\neq 0$.
\item $\left[V_{\lambda_1}\otimes V_{\lambda_2}\otimes V_{\lambda_3}\otimes{\det}^{-n}\right]^{GL_m}\neq 0$.
\end{itemize}
\end{coro}

\subsection{Duality II }\label{sec:duality-2}

Taking the transpose of the Young diagram defines a bijective map 
$\lambda \subset m\times n \longmapsto \lambda^\vee  \subset n\times m$. The following lemma is useful in our computations.

\begin{lem}
If the partition $\lambda\subset m\times n $ is equal to $\lambda(I)$ then $\lambda^\vee=\lambda(I^\vee)$
where $I^\vee=\widetilde{(I^c)}$.
\end{lem}

The canonical bilinear form on $\Cbb^{n+m}$ permits to define the map
$\delta: \Gbb_{n,m}\longrightarrow \Gbb_{m,n}$ that sends $F$ to $F^\perp$. 
Let $\delta^*: H^*(\Gbb_{m,n})\to H^*(\Gbb_{n,m})$ denotes the pullback map in cohomology.

\begin{lem}\label{lem:tilde-cohomology}
For any partition $\lambda\subset m\times n$, we have $\delta^*(\sigma_\lambda)=\sigma_{\lambda^\vee}$.
\end{lem}

\subsection{Morphism $\phi_{m,n}$ : geometric definition}

To a polynomial representation $V\in R_+(GL_m)$, we associated the polynomial map 
$\pi_V: M_{m,m}\longrightarrow \End(V)$ and the invariant polynomial $\chi_V(X):={\rm Tr}_V\left(\pi_V(X)\right)$.

Let $\Ebb_{m,n}\to \Gbb_{m,n}$ be the canonical vector bundle of rank $m$. Let $\Omega_{m,n}\in \Acal^2(\Gbb_{m,n},\End(\Ebb_{m,n}))$ be its curvature. 
The Chern-Weil homomorphism associates to the invariant polynomial $\chi_V$ the closed form 
$\chi_V\left(\frac{i}{2\pi}\Omega_{m,n}\right)$ of even degree on $\Gbb_{m,n}$. We denote by $H^*(\Gbb_{m,n})$ the de Rham cohomology of $\Gbb_{m,n}$. 
We have a natural (injective) morphism $H^*(\Gbb_{m,n},\Zbb)\to H^*(\Gbb_{m,n})$.

Here is a geometric definition of the map $\phi_{m,n}$ \cite{Tamvakis04}.

\begin{theo}\label{theo:phi-second-definition}
For any $V\in R_+(GL_m)$, $\phi_{m,n}(V)\in H^*(\Gbb_{m,n})$ is the class defined by the closed form $\chi_V\left(\frac{i}{2\pi}\Omega_{m,n}\right)$.
\end{theo}

If $\Vcal\to N$ is a complex vector bundle, we denote by $c_k(\Vcal)$ its $k$-Chern class. In the next lemma, we recall the computation of the Chern classes of the vector bundles 
$\Ebb_{m,n}$ and $\Ebb^\perp_{m,n}$.

\begin{lem} The following relations holds in  $H^*(\Gbb_{m,n})$.
\begin{align*}
c_k(\Ebb_{m,n}) &=
\begin{cases}
\sigma_{1^k}\hspace{6mm} {\rm if}\quad 1\leq k\leq m,\\
0\hspace{1cm}{\rm if}\quad k>m.
\end{cases}
\\
c_k(\Ebb_{m,n}^\perp) &=
\begin{cases}
\sigma_{k}\hspace{8mm} {\rm if}\quad 1\leq k\leq n,\\
0\hspace{1cm}{\rm if}\quad k>n.
\end{cases}
\end{align*}
\end{lem}

{\em Proof :} If $k>m={\rm rank}(\Ebb_{m,n})$, then $c_k(\Ebb_{m,n})=0$. If $1\leq k\leq m$, then 
$c_k(\Ebb_{m,n})=\phi_{m,n}({\bigwedge}^k\Cbb^m) =\sigma_{1^k}$. For the second point, let us use the isomorphism 
$\delta:\Gbb_{m,n}\to \Gbb_{n,m}$. We see that the vector bundle $\Ebb_{m,n}^\perp$ is isomorphic to 
$\delta^{-1}(\Ebb_{n,m})$. Then $c_k(\Ebb_{m,n}^\perp)=\delta^*(c_k(\Ebb_{n,m}))=\delta^*(\sigma_{1^k})=\sigma_k$  
for any $1\leq k\leq n$. $\Box$

\subsection{Cauchy formula}

We fix some integers $m,n,m',n'\geq 1$.

We consider the vector bundles  $\Ebb_{m,n}\to \Gbb_{m,n}$ and $\Ebb_{m',n'}\to \Gbb_{m',n'}$. We can form the bundles
$\Ebb_{m,n}\boxtimes\Ebb_{m',n'}$ and $\Ebb_{m,n}^\perp \boxtimes\Ebb_{m',n'}^\perp$ on $\Gbb_{m,n}\times\Gbb_{m',n'}$. 
The purpose of this section is the computation of their Euler classes.

For any partition $\lambda\subset m'\times m$, we define $\tilde{\lambda}=\widehat{\lambda}^\vee\subset m\times m'$.

\begin{prop} The following relation holds in $H^{2mm'}(\Gbb_{m,n}\times\Gbb_{m',n'})$ :
$$
\Eul(\Ebb_{m,n}\boxtimes\Ebb_{m',n'})=\sum_{\lambda\,\subset\, m'\times m}\sigma_{\tilde{\lambda}}\otimes \sigma_{\lambda}, 
$$
where $\sigma_\lambda\in H^{2|\lambda|}(\Gbb_{m',n'})$ and $\sigma_{\tilde{\lambda}}\in H^{2(mm'-|\lambda|)}(\Gbb_{m,n})$. For a partition 
$\lambda\subset m'\times m$, the product $\sigma_{\tilde{\lambda}}\otimes \sigma_{\lambda}$ does not vanish if only if the following conditions hold :
\begin{itemize}
\item $n'\geq \lambda_1$,
\item $\sharp\{1\leq k\leq m',\lambda_k=m\}\geq m'-n$.
\end{itemize}
\end{prop}

\medskip

{\em Proof : }
The Euler class $\Eul(\Ebb_{m,n}\boxtimes\Ebb_{m',n'})$ is equal to the top Chern class 
$c_{mm'}(\Ebb_{m,n}\boxtimes\Ebb_{m',n'})\in H^{2mm'}(\Gbb_{m,n}\times\Gbb_{m',n'})$. The curvature of the vector bundle 
$\Ebb_{m,n}\,\boxtimes\,\Ebb_{m',n'}$ is equal to $\Omega_{m,n}\,\otimes\, {\rm Id}' + {\rm Id}\,\otimes\, \Omega_{m',n'}$ where 
${\rm Id}\in \End(\Ebb_{m,n})$ and ${\rm Id}'\in \End(\Ebb_{m',n'})$ are the identity maps. In order to compute 
$$
c_{mm'}(\Ebb_{m,n}\boxtimes\Ebb_{m',n'})=\det\left(\tfrac{i}{2\pi}\Omega_{m,n}\,\otimes\, {\rm Id}' + {\rm Id}\,\otimes\, \tfrac{i}{2\pi}\Omega_{m',n'}\right)
$$
we use the following Cauchy formula (see \cite{Macdonal92}) 
$$
\prod_{\stackrel{1\leq i\leq m}{1\leq j\leq m'}}(x_i+x'_j)=\sum_{\lambda\subset m'\times m}s_{\tilde{\lambda}}(x) s_\lambda(x').
$$
The previous relation implies that 
$$
 \det\left(X\,\otimes\, {\rm Id}' + {\rm Id}\,\otimes\, X'\right)=\!\sum_{\lambda\,\subset \,m'\times m}\chi_{V_{\tilde{\lambda}}}(X) \chi_{V_\lambda}(X'),
$$
for all $(X,X')\in M_{m,m}\times M_{m',m'}$. Finally, we obtain thank to Theorem \ref{theo:phi-second-definition}, the following relation 
\begin{align*}
c_{mm'}(\Ebb_{m,n}\boxtimes\Ebb_{m',n'}) &= 
\sum_{\lambda\,\subset\, m'\times m} \chi_{V_{\tilde{\lambda}}}(\tfrac{i}{2\pi}\Omega_{m,n})\chi_{V_\lambda}(\tfrac{i}{2\pi}\Omega_{m',n'})\\
&= \sum_{\lambda\,\subset\, m'\times m}\sigma_{\tilde{\lambda}}\otimes\sigma_\lambda.
\end{align*}

Let us analyse when $\sigma_{\tilde{\lambda}}\otimes\sigma_\lambda\neq 0$. From the definition, we see that for any partition 
$\lambda\subset m'\times m$, we have $\tilde{\lambda}_j=\sharp\{1\leq k\leq m', \lambda_k\leq m-j\}$, $\forall j\in\{1,\ldots, m\}$.
In particular we get $\tilde{\lambda}_1=m'-\sharp\{1\leq k\leq m', \lambda_k =m\}$ since $\lambda_k\leq m$, $\forall k$.

The element $\sigma_\lambda\in H^{2|\lambda|}(\Gbb_{m',n'})$ does not vanish if and only if $\lambda_1\leq n'$, and 
$\sigma_{\tilde{\lambda}}\in H^{2(mm'-|\lambda|)}(\Gbb_{m,n})$ does not vanish if and only if $\tilde{\lambda}_1\leq n$. 
Our proof is completed since the relation $\tilde{\lambda}_1\leq n$ is equivalent to  $\sharp\{1\leq k\leq m',\lambda_k=m\}\geq m'-n$. 
$\Box$

\medskip

We apply our formula to the case where $p\geq q\geq 1$, $1< r< p$ and $1< s< q$. 
The following relation holds in $H^{2(p-r)s}(\Gbb_{p-r,r}\times\Gbb_{s,q-s})$ :
\begin{equation}\label{eq:euler}
\Eul(\Ebb_{p-r,r}\boxtimes\Ebb_{s,q-s})=\sum_{\lambda\,\subset\, s\times p-r}\sigma_{\tilde{\lambda}}\otimes \sigma_{\lambda}, 
\end{equation}
where $\sigma_\lambda\in H^{2|\lambda|}(\Gbb_{p-r,r})$ and $\sigma_{\tilde{\lambda}}\in H^{2((p-r)s-|\lambda|)}(\Gbb_{s,q-s})$. For a partition 
$\lambda\subset s\times p-r$, the product $\sigma_{\tilde{\lambda}}\otimes \sigma_{\lambda}$ does not vanish if only if $q-s\geq \lambda_1$,
and $\sharp\{1\leq k\leq s,\lambda_k=p-r\}\geq s-r$.

\begin{lem}\label{lem:euler-0}
Suppose that $p\geq q$. Then $\Eul(\Ebb_{p-r,r}\boxtimes\Ebb_{s,q-s})=0$ if $s>r$.
\end{lem}

{\em Proof :} Suppose that $s-r>0$. If $\Eul(\Ebb_{p-r,r}\boxtimes\Ebb_{s,q-s})\neq 0$, there exists a partition $\lambda\subset s\times p-r$ such that 
$\sigma_{\tilde{\lambda}}\otimes \sigma_{\lambda}\neq 0$: hence $q-s\geq \lambda_1$ and $\sharp\{1\leq k\leq s,\lambda_k=p-r\}\geq s-r>0$. 
We obtain $q-s\geq \lambda_1=p-r$, so $q\geq p+s-r>p$, which is in contradiction with our hypothesis. $\Box$

\bigskip

We consider now the vector bundle $\Vcal^r_s= \Ebb_{r,p-r}^\perp\boxtimes\Ebb_{q-s,s}^\perp$ on 
$\Gbb_{r,p-r}\times\Gbb_{q-s,s}$. 

\begin{prop} \label{prop:calcul-euler} Let $p\geq q\geq 1$, $1< r< p$ and $1< s< q$. 

\begin{itemize}
\item If $s>r$, then $\Eul(\Vcal^r_s)=0$
\item If $s\leq r$, the following relation holds in $H^{2(p-r)s}(\Gbb_{r,p-r}\times\Gbb_{q-s,s})$ :
$$
\Eul(\Vcal^r_s)=\sum_{\lambda\,\subset\, s\times p-r} \sigma_{\widehat{\lambda}}\otimes \sigma_{\lambda^\vee}, 
$$
where $\sigma_{\lambda^\vee}\in H^{2|\lambda|}(\Gbb_{q-s,s})$ and $\sigma_{\widehat{\lambda}}\in H^{2((p-r)s-|\lambda|)}(\Gbb_{r,p-r})$. For a partition 
$\lambda\subset s\times p-r$, the product $\sigma_{\widehat{\lambda}}\otimes \sigma_{\lambda^\vee}$ does not vanish if and only if $q-s\geq \lambda_1$.
\end{itemize}
\end{prop}

{\em Proof :} Let $\delta\times \delta : \Gbb_{r,p-r}\times\Gbb_{q-s,s}\to \Gbb_{p-r,r}\times\Gbb_{s,q-s}$ be the product of ``duality'' maps 
(see \S \ref{sec:duality-2}). Since $\Ebb_{r,p-r}^\perp\boxtimes\Ebb_{q-s,s}^\perp\simeq(\delta\times \delta)^{-1}(\Ebb_{p-r,r}\boxtimes\Ebb_{q-s,s})$, we have
$\Eul(\Vcal^r_s)=\delta^*\times \delta^*(\Eul(\Ebb_{p-r,r}\boxtimes\Ebb_{s,q-s}))$. We can use 
Lemmas \ref{lem:euler-0} and \ref{lem:tilde-cohomology}, and  (\ref{eq:euler}) to complete the proof. $\Box$

\section{Convex cone ${\rm S}(p,q)$ : equations of the facets }\label{sec:s-p-q-facets}

Theorem \ref{theo:delta-cas-general} tells us that an element $(A,B,C)\in(\Ccal_p\times\Ccal_q)^3$ belongs to ${\rm S}(p,q)=\Delta((GL_p\times GL_q)^2\times \Cbb^p\otimes\Cbb^q)$ if and only  
$$
\boxed{|A'| + |B'| + |C'|=|A''| + | B'' | + | C'' |}
$$ 
and 
\begin{equation}\label{eq:equation-facet-bis}
\boxed{|\,A'\,|_{I'} + |\,B'\,|_{J'} + |\,C'\,|_{K'}\leq |\,A''\,|_{I''} + |\,B''\,|_{J''} + |\,C''\,|_{K''}},
\end{equation}
for any couple $(r,s)\in \{0,\ldots,p\}\times \{0,\ldots,q\}-\{(p,q), (0,0)\}$, for any $I',J',K'\in\Pcal^p_r$ and any $I'',J'',K''\in\Pcal^p_s$, such that the product 
\begin{equation}\label{eq:cohomological-condition}
\left(\sigma_{\lambda(I')}\otimes \sigma_{\lambda((I'')^c)}\right)\cdot \left(\sigma_{\lambda(J')}\otimes \sigma_{\lambda((J'')^c)}\right)\cdot 
\left(\sigma_{\lambda(K')}\otimes \sigma_{\lambda((K'')^c)}\right)\cdot {\rm Eul}(\Vcal^r_s)
\end{equation}
 is a non zero multiple of $[pt]\in H^{\max}\left(\Gbb(r,p-r)\times \Gbb(q-s,s),\Zbb\right)$.

In the following sections, we study each case according to the parameter $(r,s)$.

\subsection{$r=0$ and $s=q$}

Here, $\Gbb(r,p-r)\times\Gbb(q-s,s)=\{pt\}$ and $\Vcal^0_q\neq 0$, hence the cohomological condition (\ref{eq:cohomological-condition}) does not hold. 

\subsection{$r=p$ and $s=0$}

Here, $\Gbb(r,p-r)\times\Gbb(q-s,s)=\{pt\}$ and $\Vcal^p_0= 0$, hence the cohomological condition (\ref{eq:cohomological-condition}) holds. Relation 
(\ref{eq:equation-facet-bis}) becomes
$\boxed{|\,A'\,| + |\,B'\,| + |\,C'\,|\leq 0}$.

\subsection{$0<r<p$ and $s=0$}
 
 Here $\Vcal^r_0=0$ and the cohomological condition (\ref{eq:cohomological-condition})  becomes 
 $\sigma_{\lambda(I')}\cdot \sigma_{\lambda(J')}\cdot \sigma_{\lambda(K')}=k [pt],\ k\geq 1$, in $H^*(\Gbb(r,p-r))$ which is equivalent to asking that 
 $$
 \left[V_{\lambda(I')}\otimes V_{\lambda(J')}\otimes V_{\lambda(K')}\otimes {\det}^{-(p-r)}\right]^{GL_r}\neq 0.
 $$

Here inequality (\ref{eq:equation-facet-bis}) becomes $\boxed{|\,A'\,|_{I'} + |\,B'\,|_{J'} + |\,C'\,|_{K'}\leq 0}$.

\subsection{$0<r<p$ and $s=q$}
 Here $\Vcal^r_q$ is the vector bundle $\Ebb_{r,p-r}^\perp\otimes \Cbb^q$ on $\Gbb(r,p-r)$. Since 
 $$
 \Eul(\Vcal^r_q)=\Eul(\Ebb_{r,p-r}^\perp)^q=\left(c_{p-r}(\Ebb_{r,p-r}^\perp)\right)^q= (\sigma_{p-r})^q\quad \in H^*(\Gbb(r,p-r)),
 $$
condition (\ref{eq:cohomological-condition}) is 
$\sigma_{\lambda(I')}\cdot \sigma_{\lambda(J')}\cdot \sigma_{\lambda(K')}\cdot(\sigma_{p-r})^q=k [pt],\ k\geq 1$. If we take the image of the previous relation 
through $\delta: H^*(\Gbb(r,p-r))\to H^*(\Gbb(p-r,r))$, we obtain
$$
\sigma_{\lambda(\widetilde{(I')^c})}\cdot \sigma_{\lambda(\widetilde{(J')^c})}\cdot \sigma_{\lambda(\widetilde{(K')^c})}\cdot(\sigma_{1^{p-r}})^q=k [pt],\ k\geq 1,\quad {\rm in}\quad H^*(\Gbb(p-r, r))
$$
that is equivalent to 
$$
 \left[V_{\lambda(\widetilde{(I')^c})}\otimes V_{\lambda(\widetilde{(J')^c})}\otimes V_{\lambda(\widetilde{(K')^c})}\otimes {\det}^{q-r}\right]^{GL_{p-r}}\neq 0.
$$

Here inequality (\ref{eq:equation-facet-bis}) is equivalent to $\boxed{|\,A'\,|_{(I')^c} + |\,B'\,|_{(J')^c} + |\,C'\,|_{(K')^c}\geq 0}$.

\subsection{$r=0$ and $0<s<q$}
 Here $\Vcal^0_s$ is the vector bundle $\Cbb^p\otimes \Ebb_{q-s,s}^\perp$ on $\Gbb(p-s,s)$. Since 
 $\Eul(\Vcal^0_s)=(\sigma_{s})^p$, condition (\ref{eq:cohomological-condition}) is 
$\sigma_{\lambda((I'')^c)}\cdot \sigma_{\lambda((J'')^c)}\cdot \sigma_{\lambda((K'')^c)}\cdot(\sigma_{s})^p=k [pt],\ k\geq 1$. If we take the image of the previous relation 
through $\delta: H^*(\Gbb(q-s,s))\to H^*(\Gbb(s,q-s))$, we obtain
$$
\sigma_{\lambda(\widetilde{I''})}\cdot \sigma_{\lambda(\widetilde{J''})}\cdot \sigma_{\lambda(\widetilde{K''})}\cdot(\sigma_{1^{s}})^p=k [pt],\ k\geq 1,\quad {\rm in}\quad H^*(\Gbb(s,q-s))
$$
that is equivalent to  
$$
 \left[V_{\lambda(\widetilde{I''})}\otimes V_{\lambda(\widetilde{J''})}\otimes V_{\lambda(\widetilde{K''})}\otimes {\det}^{p-(q-s)}\right]^{GL_{s}}\neq 0.
$$

Here inequality (\ref{eq:equation-facet-bis}) becomes $\boxed{|\,A''\,|_{I''} + |\,B'\,|_{J''} + |\,C'\,|_{K''}\geq 0}$.

\subsection{$r=p$ and $0<s<q$}

 Here $\Vcal^p_s=0$ and the cohomological condition (\ref{eq:cohomological-condition})  becomes 
 $\sigma_{\lambda((I'')^c)}\cdot \sigma_{\lambda((J'')^c)}\cdot \sigma_{\lambda((K'')^c)}=k [pt],\ k\geq 1$, in $H^*(\Gbb(q-s,s))$ that is equivalent to 
 $$
 \left[V_{\lambda((I'')^c)}\otimes V_{\lambda((J'')^c)}\otimes V_{\lambda((K'')^c)}\otimes {\det}^{-s}\right]^{GL_{q-s}}\neq 0.
 $$
 
Here inequality (\ref{eq:equation-facet-bis}) is  equivalent to $\boxed{|\,A''\,|_{(I'')^c} + |\,B''\,|_{(J'')^c} + |\,C''\,|_{(K'')^c}\leq 0}$.

\subsection{$0<r<p$ and $0<s<q$}
 
We know that $\Eul(\Vcal^r_s)=0$ if $s>r$, hence condition (\ref{eq:cohomological-condition})  does not hold if $s>r$. 
Suppose now that $s\leq r$. Thanks to Proposition \ref{prop:calcul-euler}, we see that condition (\ref{eq:cohomological-condition}) is equivalent 
to asking that the sum  
$$
\sum_{\lambda\subset s\times p-r}\left(\sigma_{\lambda((I')^c)}\cdot \sigma_{\lambda((J')^c)}\cdot \sigma_{\lambda((K')^c)}\cdot\sigma_{\widehat{\lambda}}\right)
\otimes \left(\sigma_{\lambda(I'')}\cdot \sigma_{\lambda(J'')}\cdot \sigma_{\lambda(K'')}\cdot\sigma_{\lambda^\vee}\right)
$$
is equal to $k[pt],k\geq 1$ in $H^*(\Gbb(r,p-r))\otimes H^*(\Gbb(q-s,s))$. Hence (\ref{eq:cohomological-condition}) holds if and only if there exists 
a partition $\lambda\subset s\times p-r$ such that 
\begin{align}
\sigma_{\lambda(I')}\cdot \sigma_{\lambda(J')}\cdot \sigma_{\lambda(K')}\cdot\sigma_{\widehat{\lambda}} &= \ell'[pt]  
& {\rm in}\  H^*(\Gbb(r,p-r),\Zbb),\label{eq:coho-1}
\\
\sigma_{\lambda((I'')^c)}\cdot \sigma_{\lambda((J'')^c)}\cdot \sigma_{\lambda((K'')^c)}\cdot
\sigma_{\lambda^\vee}& = \ell'' [pt] & {\rm in} \ H^*(\Gbb(q-s,s),\Zbb),\label{eq:coho-2}
\end{align}
for some $\ell',\ell''\geq 1$.
\begin{lem}
If $s\leq r$, then (\ref{eq:coho-1}) and (\ref{eq:coho-2}) hold if and only if there exists a partition $\mu\subset s\times p-r$, such that 
both conditions hold
 \begin{align}
 \left[V_{\lambda(I')}\otimes V_{\lambda(J')}\otimes V_{\lambda(K')}\otimes V_\mu\otimes {\det}^{-(p-r)}\right]^{GL_r}& \neq 0, \label{eq:tensor-1}
\\
 \left[V_{\lambda(I'')}\otimes V_{\lambda(J'')}\otimes V_{\lambda(K'')}\otimes V_\mu\otimes {\det}^{-(p-r)-2(q-s)}\right]^{GL_s}& \neq 0.\label{eq:tensor-2}
\end{align}
\end{lem}

{\em Proof :}  Let $\lambda\subset s\times p-r$ satisfying (\ref{eq:coho-1}) and (\ref{eq:coho-2}). Let $\mu=\widehat{\lambda}\subset s\times p-r$. 
Condition (\ref{eq:coho-1}) is then equivalent\footnote{Since the length of $\mu$ is less than $r\geq s$, $V_\mu$ is well defined as an 
irreducible representation of $GL_r$.} to (\ref{eq:tensor-1}).

Take now the image of (\ref{eq:coho-2}) through the map $\delta: H^*(\Gbb(q-s,s),\Zbb)\to H^*(\Gbb(s,q-s),\Zbb)$ : we obtain
$\sigma_{\lambda(\widetilde{I''})}\cdot \sigma_{\lambda(\widetilde{J''})}\cdot \sigma_{\lambda(\widetilde{K''})}\cdot
\sigma_{\lambda}= \ell'' [pt] $, $\ell''\geq 1$  in $H^*(\Gbb(s,q-s),\Zbb)$, that is equivalent to 
\begin{equation}\label{eq:intermediaire}
[V_{\lambda(\widetilde{I''})}\otimes V_{\lambda(\widetilde{J''})}\otimes V_{\lambda(\widetilde{K''})}\otimes V_{\lambda}\otimes {\det}^{-(q-s)}]^{GL_s}\neq 0.
\end{equation}
Since we have the following relations between representations of $GL_s$ : 
\begin{enumerate}
\item $V_{\lambda(\widetilde{I''})}=(V_{\lambda(I'')})^*\otimes {\det}^{q-s}$ (with the same relations for $J''$ and $K''$),
\item $V_{\lambda}=(V_\mu)^*\otimes  {\det}^{p-r}$,
\end{enumerate}
 condition (\ref{eq:intermediaire}) is equivalent to (\ref{eq:tensor-2}). $\Box$

\subsection{Summary}\label{sec:summary}

Let us summarize the computations done in the previous sections. 

An element $(A,B,C)\in(\Ccal_p\times\Ccal_q)^3$ belongs to ${\rm S}(p,q)$ if and only 
the following conditions hold
\begin{itemize}
\item $\boxed{|A'| + |B'| + |C'|=|A''| + | B'' | + | C'' |}$.
\item $\boxed{|\,A'\,| + |\,B'\,| + |\,C'\,|\leq 0}$.
\item For any $0<r<p$, for any $I',J',K'\in \Pcal^p_r$, we have :
\begin{align*}
\boxed{|\,A'\,|_{I'} + |\,B'\,|_{J'} + |\,C'\,|_{K'}\leq 0} \quad &{\rm if} \  
\left[V_{\lambda(I')}\otimes V_{\lambda(J')}\otimes V_{\lambda(K')}\otimes {\det}^{-p+r}\right]^{GL_r}\neq 0, \\
\boxed{|\,A'\,|_{I'} + |\,B'\,|_{J'} + |\,C'\,|_{K'}\geq 0} \quad & {\rm if}\  \left[V_{\lambda(\widetilde{I'})}\otimes 
V_{\lambda(\widetilde{J'})}\otimes V_{\lambda(\widetilde{K'})}\otimes {\det}^{q-p+r}\right]^{GL_{r}}\neq 0.
\end{align*}
\item For any $0<s<q$, for any $I'',J'',K''\in \Pcal^q_s$, we have :
\begin{align*}
\boxed{|\,A''\,|_{I''} + |\,B''\,|_{J''} + |\,C''\,|_{K''}\leq 0} \quad & {\rm if}\  \left[V_{\lambda(I'')}\otimes V_{\lambda(J'')}\otimes 
V_{\lambda(K'')}\otimes {\det}^{-q+s}\right]^{GL_{s}}\neq 0.\\
\boxed{|\,A''\,|_{I''} + |\,B''\,|_{J''} + |\,C''\,|_{K''}\geq 0} \quad &{\rm if} \  \left[V_{\lambda(\widetilde{I''})}\otimes V_{\lambda(\widetilde{J''})}
\otimes V_{\lambda(\widetilde{K''})}\otimes {\det}^{p-q+s}\right]^{GL_{s}}\neq 0.
\end{align*}

\item For any $(r,s)\in [p-1]\times [q-1]$ with $r\geq s$, for any $I,J,K\in \Pcal^p_r\times\Pcal^q_s$, we have 
$$
\boxed{|\,A'\,|_{I'} + |\,B'\,|_{J'} + |\,C'\,|_{K'}\leq |\,A''\,|_{I''} + |\,B''\,|_{J''} + |\,C''\,|_{K''}}
$$
if there exists a partition $\mu\subset s\times p-r$, such that 
both conditions hold
 \begin{align*}
 \left[V_{\lambda(I')}\otimes V_{\lambda(J')}\otimes V_{\lambda(K')}\otimes V_\mu\otimes {\det}^{-(p-r)}\right]^{GL_r}& \neq 0, 
\\
 \left[V_{\lambda(I'')}\otimes V_{\lambda(J'')}\otimes V_{\lambda(K'')}\otimes V_\mu\otimes {\det}^{-(p-r)-2(q-s)}\right]^{GL_s}& \neq 0.
\end{align*}
\end{itemize}

\subsection{Proof of the main result}

In the previous section, we described the facets of the convex cone ${\rm S}(p,q)$. We will now exploit the fact that 
$(A,B,C)\in\horn(p,q)$ if and only if 
$$
\Theta(A,B,C)=\left(\begin{pmatrix}
A' \\
(A')^*
\end{pmatrix},
\begin{pmatrix}
B' \\
(B'')^*
\end{pmatrix},
\begin{pmatrix}
(C')^* \\
C''
\end{pmatrix}
\right)\in {\rm S}(p,q).
$$

In what follows, we crucially use that the semigroups $\horn^\Zbb(n)$ and $\horn^{\Zbb}(p,q)$ are saturated. Hence 
$\horn^{\Zbb}(n)=\horn(n)\cap (\wedge^+_n)^3$ and  \break $\horn^{\Zbb}(p,q)=\horn(p,q)\cap (\wedge^+_p\times\wedge^+_q)^ 3$. 
Recall also the identity $|A^*|_{I}=-|A|_{\widetilde{I}}$ which will be used several times. 

Let us compute the image of the facets of ${\rm S}(p,q)$ through the linear map $\Theta$ :
\begin{itemize}
\item The image of $|A'| + |B'| + |C'|=|A''| + | B'' | + | C'' |$ through $\Theta$ is $\boxed{|\,A\,| + |\,B\,| = |\,C\,|}$.
\item The image of the half space $|\,A'\,| + |\,B'\,| + |\,C'\,|\leq 0$ through the map $\Theta$ is $\boxed{|\,A'\,| + |\,B'\,| \leq |\,C'\,|}$.
\item The image of the half space  $|\,A'\,|_{I'} + |\,B'\,|_{J'} + |\,C'\,|_{K'}\leq 0$ through the map $\Theta$ is 
$\boxed{|\,A'\,|_{I'} + |\,B'\,|_{J'} \leq |\,C'\,|_{\widetilde{K'}}}$, 
and condition 
$$
\left[V_{\lambda(I')}\otimes V_{\lambda(J')}\otimes V_{\lambda(K')}\otimes {\det}^{-p+r}\right]^{GL_r}\neq 0
$$
 is equivalent to $(\lambda(I'),\lambda(J'),\lambda(\widetilde{K'}))\in\horn(r)$.

\item The image of of the half space $|\,A'\,|_{I'} + |\,B'\,|_{J'} + |\,C'\,|_{K'}\geq 0$ through the map 
$\Theta$ is $\boxed{|\,A'\,|_{I'} + |\,B'\,|_{J'} \geq |\,C'\,|_{\widetilde{K'}}}$, and  condition 
$$
\left[V_{\lambda(\widetilde{I'})}\otimes 
V_{\lambda(\widetilde{J'})}\otimes V_{\lambda(\widetilde{K'})}\otimes {\det}^{q-p+r}\right]^{GL_{r}}\neq 0
$$
 is equivalent to 
$(\lambda(I'),\lambda(J'),\lambda(\widetilde{K'})+(q\!+\!p\!-\!r)\mathbf{1}_r)\in\horn(r)$.
\item The image of of the half space  $|\,A''\,|_{I''} + |\,B''\,|_{J''} + |\,C''\,|_{K''}\leq 0$  through the map 
$\Theta$ is $\boxed{|\,A''\,|_{\widetilde{I''}} + |\,B''\,|_{\widetilde{J''}} \geq |\,C''\,|_{K''}}$ and condition 
$$
\left[V_{\lambda(I'')}\otimes V_{\lambda(J'')}\otimes 
V_{\lambda(K'')}\otimes {\det}^{-q+s}\right]^{GL_{s}}\neq 0
$$
is equivalent to $(\lambda(\widetilde{I''}),\lambda(\widetilde{J''}),\lambda(K'')+(q\!-\!s)\mathbf{1}_s)\in\horn(s)$.
\item The image of of the half space  $|\,A''\,|_{I''} + |\,B''\,|_{J''} + |\,C''\,|_{K''}\geq 0$  through the map 
$\Theta$ is $\boxed{|\,A''\,|_{\widetilde{I''}} + |\,B''\,|_{\widetilde{J''}} \leq |\,C''\,|_{K''}}$ and condition 
$$
\left[V_{\lambda(\widetilde{I''})}\otimes V_{\lambda(\widetilde{J''})}
\otimes V_{\lambda(\widetilde{K''})}\otimes {\det}^{p-q+s}\right]^{GL_{s}}\neq 0
$$
is equivalent to $(\lambda(\widetilde{I''}),\lambda(\widetilde{J''}),\lambda(K'')-p\mathbf{1}_s)\in\horn(s)$.

\item The image of of the half space 
$$
|\,A'\,|_{I'} + |\,B'\,|_{J'} + |\,C'\,|_{K'}\leq |\,A''\,|_{I''} + |\,B''\,|_{J''} + |\,C''\,|_{K''}
$$
through the map $\Theta$ is 
$$
\boxed{|\,A'\,|_{I'} + |\,A''\,|_{\widetilde{I''}}+ |\,B'\,|_{J'} + |\,B''\,|_{\widetilde{J''}} \leq  |\,C'\,|_{\widetilde{K'}} + |\,C''\,|_{K''}} 
$$
and conditions 
 \begin{align*}
 \left[V_{\lambda(I')}\otimes V_{\lambda(J')}\otimes V_{\lambda(K')}\otimes V_\mu\otimes {\det}^{-(p-r)}\right]^{GL_r}& \neq 0, 
\\
 \left[V_{\lambda(I'')}\otimes V_{\lambda(J'')}\otimes V_{\lambda(K'')}\otimes V_\mu\otimes {\det}^{-(p-r)-2(q-s)}\right]^{GL_s}& \neq 0.
\end{align*}
are equivalent to 
 \begin{align*}
 \left[V_{\lambda(I')}\otimes V_{\lambda(J')}\otimes V_{\lambda(\widetilde{K'})}^*\otimes V_\mu\right]^{GL_r}& \neq 0, 
\\
 \left[V_{\lambda(\widetilde{I''})}\otimes V_{\lambda(\widetilde{J''})}\otimes V_{\lambda(K'')}^*\otimes {\det}^{(p-r)}\otimes V_\mu^*\right]^{GL_s}& \neq 0.
\end{align*}
The existence of a partition $\mu\in s\times p-r$ satisfying the previous relations  is equivalent to asking that 
$$
\left(\begin{pmatrix}
\lambda(I') \\
\lambda(\widetilde{I''})
\end{pmatrix},
\begin{pmatrix}
\lambda(J')\\
\lambda(\widetilde{J''})
\end{pmatrix},
\begin{pmatrix}
\lambda(\widetilde{K'}) \\
\lambda(K'')-(p-r)\mathbf{1}_s
\end{pmatrix}
\right)\in (\wedge^+_r\times \wedge^+_s)^3
$$
belongs to $\horn(r,s)$.
\end{itemize}

The proof of Theorem \ref{theo:main} is complete.


\end{document}